\documentclass[a4paper,10pt]{scrartcl}

\usepackage{orcidlink}
\usepackage{comment}
\usepackage [autostyle, english = british]{csquotes}
\MakeOuterQuote{"}
\usepackage{wrapfig}
\usepackage[english]{babel}
\usepackage{blindtext}

\usepackage{amssymb,bbm}

\usepackage{amsmath,amsthm}
\newtheorem{theorem}{Theorem}
\newtheorem{lemma}{Lemma}

\newtheorem{remark}{Remark}
\newtheorem{proposition}{Proposition}

\usepackage{cleveref,caption,subcaption,hyperref}

\usepackage{algorithm} 
\usepackage{algcompatible}
\Crefname{ALC@unique}{Line}{Lines}

\usepackage{graphicx}
\usepackage{xcolor}
\usepackage{pgfplots}
\pgfplotsset{compat=newest, width=0.5\columnwidth} 

\newlength\fwidth
\newlength\fheight
\tikzset{
    dot/.style 2 args={fill, circle, inner sep=1pt, label={#1:\scriptsize #2}}
}

\definecolor{cPOD}{HTML}{D7191C}  
\definecolor{cOI3}{HTML}{253494}  
\definecolor{cOI2}{HTML}{2C7FB8}  
\definecolor{cOI1}{HTML}{41B6C4}  

\newenvironment{keywords}%
   {\begin{trivlist}\item[]{\bfseries\sffamily Keywords:}\ }
   {\end{trivlist}}

\newenvironment{MSCcodes}%
   {\begin{trivlist}\item[]{\bfseries\sffamily Mathematics subject classification:}\ }
   {\end{trivlist}}

\newcommand{\email}[1]{{\normalfont Email:~\texttt{\href{mailto:#1}{#1}}}}

\DeclareMathOperator*{\Cov}{\mathrm{Cov}}
\DeclareMathOperator*{\Var}{\mathrm{Var}}

\DeclareMathOperator*{\argmin}{arg\,min}
\DeclareMathOperator*{\diag}{\textbf{diag}}
\newcommand{\R}{\ensuremath{\mathbb{R}}}
\newcommand{\Rn}{\ensuremath{\mathbb{R}^n}}

\newcommand{\Rnn}{\ensuremath{\mathbb{R}^{n\times n}}}
\newcommand{\Rnm}{\ensuremath{\mathbb{R}^{n\times m}}}

\newcommand*{\rmd}{\mathop{}\!\mathrm{d}}
\newcommand*{\rmE}{\ensuremath{\mathrm{E}}}
\newcommand*{\E}[1]{\mathop{}\!\mathbb{E}\left[ #1 \right]}

\newcommand*{\bbX}{\mathop{}\!\mathbb{X}}

\newcommand{\A}{\ensuremath{\mathrm{A}}}
\newcommand{\B}{\ensuremath{\mathrm{B}}}
\newcommand{\I}{\ensuremath{\mathrm{I}}}
\newcommand{\N}{\ensuremath{\mathrm{N}}}
\renewcommand{\H}{\ensuremath{\mathrm{H}}}
\newcommand{\M}{\ensuremath{\mathrm{M}}}
\renewcommand{\O}{\ensuremath{\mathrm{O}}}
\newcommand{\K}{\ensuremath{\mathrm{K}}}
\renewcommand{\S}{\ensuremath{\mathrm{S}}}
\newcommand{\D}{\ensuremath{\mathrm{D}}}
\newcommand{\F}{\ensuremath{\mathrm{F}}}

\newcommand{\U}{\ensuremath{\mathrm{U}}}
\newcommand{\T}{\ensuremath{\mathrm{T}}}

\newcommand{\V}{\ensuremath{\mathrm{V}}}
\newcommand{\C}{\ensuremath{\mathrm{C}}}
\newcommand{\W}{\ensuremath{\mathrm{W}}}

\newcommand{\Z}{\ensuremath{\mathrm{Z}}}

\newcommand{\rmR}{\ensuremath{\mathrm{R}}}
\newcommand{\rmP}{\ensuremath{\mathrm{P}}}
\newcommand{\hA}{\ensuremath{\hat{\mathrm{A}}}}
\newcommand{\hB}{\ensuremath{\hat{\mathrm{B}}}}
\newcommand{\hN}{\ensuremath{\hat{\mathrm{N}}}}
\newcommand{\hH}{\ensuremath{\hat{\mathrm{H}}}}
\newcommand{\hM}{\ensuremath{\hat{\mathrm{M}}}}
\newcommand{\hO}{\ensuremath{\hat{\mathrm{O}}}}
\newcommand{\hS}{\ensuremath{\hat{\mathrm{S}}}}
\newcommand{\tO}{\ensuremath{\tilde{\mathrm{O}}}}
\newcommand{\hK}{\ensuremath{\hat{\mathrm{K}}}}
\newcommand{\hPsi}{\ensuremath{\hat{\Psi}}}

\newcommand{\hX}{\ensuremath{\hat{X}}}

\newcommand*{\equal}{=}
\newcommand{\hb}[1]{\ensuremath{\hat{\bar{#1}}}}
\newcommand{\db}[1]{\ensuremath{\dot{\bar{#1}}}}

\title{Learning Stochastic Reduced Models from Data: A Nonintrusive Approach}
\author{M. A. Freitag\,\orcidlink{0000-0002-4539-2162}\thanks{Institut für Mathematik, Universität Potsdam, 14476 Potsdam, Germany (\email{melina.freitag@uni-potsdam.de}, \email{jan.martin.nicolaus@uni-potsdam.de}).} \and J. M. Nicolaus\,\orcidlink{0009-0005-1325-3626}\footnotemark[1] \thanks{Institut für Mathematik, Universität Rostock, 18051 Rostock, Germany, (\email{martin.redmann@uni-rostock.de}).} \and M. Redmann\orcidlink{0009-0000-3009-3138}\footnotemark[2]}
\date{\today}

\begin{document}

\maketitle
\begin{abstract}
A nonintrusive model order reduction method for bilinear stochastic differential equations with additive Gaussian noise is proposed.
A reduced order model (ROM) is designed in order to approximate the statistical properties of high-dimensional systems.
The drift and diffusion coefficients of the ROM are inferred from state observations by solving appropriate least-squares problems.
The closeness of the ROM obtained by the presented approach to the intrusive ROM obtained by the proper orthogonal decomposition (POD) method is investigated. 
Two generalisations of the snapshot-based dominant subspace construction to the stochastic case are presented.
Numerical experiments are provided to compare the developed approach to POD.
\end{abstract}


\section{Introduction} 
The numerical solution of (stochastic) partial differential equations (PDE) is ubiquitous throughout many fields in engineering and applied sciences. 
To obtain a numerical model, a spatial discretisation of the governing PDE, for example by finite differences, can be performed.
To achieve high accuracy and resolve small-scale phenomena, the utilised discretisation mesh has to be finely grained, resulting in a numerical model of high dimension and computationally expensive evaluations. 
In a setting where such a model has to be evaluated often, for instance, when optimising over a parameter or when a Monte-Carlo simulation is run, the computational complexity might result in a practically infeasible algorithm, due to limited computational resources or time sensitivity in an online setting.

The field of model order reduction (MOR) deals with the construction of surrogate models that are much cheaper to evaluate than the original full order model (FOM). 
In this paper, the FOM is a controlled stochastic differential equation (SDE) with bilinear drift and additive Gaussian noise. 
That is, the $n$-dimensional FOM is of the form
\begin{equation}\label{eq:SDEFOM}
\rmd X(t) = \left[\A X(t) + \B u(t) + \sum_{i=1}^m \N_{i}X(t)u_i(t)\right]\rmd t + \M \rmd W(t), \quad X(0) = X_0,
\end{equation}
where $W(t)$, $t\in[0, T]$, is a $d$-dimensional Wiener process with correlation matrix $\K\in\R^{d\times d}$ and $u:[0,T]\rightarrow\R^m$ an $m$-dimensional integrable deterministic control.
The initial condition $X_0 \in\Rn$ is assumed to be either deterministic or a Gaussian random vector.
All random variables, vectors, and processes appearing in this paper are defined on the filtered probability space \((\Omega, \mathcal{F}, \{\mathcal{F}_t\}_{ t\in[0,T]}, \mathbb{P})\).

Originating from deterministic systems ($M\equiv 0$), a wide range of reduction techniques has been developed.
Based on the accessibility of the FOM coefficients $\A$, $\N_1, \dots, \N_m \in \Rnn,\B\in\Rnm,\M\in\R^{n\times d}, \text{ and }\K\in\R^{d\times d}$, these approaches are divided into two categories.
On the one hand, intrusive model reduction requires knowledge of the system coefficients of the FOM to construct the reduced order model (ROM).
A famous projection-based method for the reduction of linear time-invariant systems from a system-theoretic perspective is balanced truncation (BT) \cite{moore_principal_1981}, which has been recently extended to the stochastic case \cite{benner2015model,Benner2015DualPO,becker2019infinite,redmann_model_2024}.
Another well-known and closely related, but data-driven, approach is the proper orthogonal decomposition (POD) method \cite{kunisch_galerkin_2001}.
Given snapshots of the system states of the FOM, the POD method constructs an $r$-dimensional subspace $\mathcal{V}_r\subset \Rn$ as the span of the $r$ leading left singular vectors $v_1,\dots,v_r$ of a snapshot matrix.
The matrix $\V_r =\begin{bmatrix}
    v_1,\dots,v_r
\end{bmatrix}\in\R^{n\times r}$ is then used to define a projection $\rmP=\V_r\V_r^T\in\Rnn$ onto $\mathcal{V}_r$.
If the FOM drift coefficients $\A,\N_1,\dots , \N_m\in\Rnn, \B\in\Rnm,$ as well as the diffusion coefficient $\M\in\R^{n\times d}$ and the correlation matrix $\K\in\R^{d\times d}$ are available, performing a Galerkin projection yields an $r$-dimensional system
\begin{equation}\label{eq:SDEROM}
\rmd X_r(t) = \left[\A_r X_r(t) + \B_r u(t) + \sum_{i=1}^m \N_{r,i}X_r(t)u_i(t)\right]\rmd t + \M_r \rmd W(t),  
\end{equation}
with $X_r(0) = X_{r,0}\in\R^r$ as the POD-ROM.
The projected coefficients and initial condition are defined by $\A_r = \V_r^T\A\V_r\in\R^{r\times r}$, $\B_r = \V_r^T \B\in\R^{r\times m}$, $\N_{r,i} = \V_r^T \N_i\V_r\in\R^{r\times r}$ and $\M_r = \V_r^T\M\in\R^{r\times d}$, as well as $X_{r,0}=\V_r^T X_0$.
The approximation to the FOM is subsequently obtained by lifting $X(t)\approx \V_rX_r(t)$, $t\in[0,T]$.
A structure preserving extension of POD to stochastic systems has recently been proposed \cite{tyranowski_data-driven_2022}.
In the presence of componentwise evaluated nonlinearities, the discrete empirical interpolation method (DEIM) \cite{chaturantabut_nonlinear_2010,chaturantabut_state_2012} provides a reduction method that can be effectively combined with POD. 
An optimisation-based intrusive MOR technique is the iterative rational Krylov algorithm (IRKA) that aims to minimise a certain error bound between the FOM and the surrogate model \cite{gugercin_mathcalh_2_2008}. It has also been generalised to SDEs, see \cite{redmann_optimization_2021}.

On the other hand, nonintrusive reduction methods construct a surrogate ROM from data without explicit knowledge of the FOM system coefficients.
Examples of nonintrusive reduction methods include the Loewner framework \cite{antoulas_model_2016,ionita_data-driven_2014,simard_nonlinear_2021,simard_construction_2024,karachalios_bilinear_2021,mayo_framework_2007} and dynamical mode decomposition \cite{schmid_dynamic_2010,williams_datadriven_2015,h_tu_dynamic_2014,kutz_dynamic_2016,rowley_spectral_2009}.
Applied to Equation \eqref{eq:SDEFOM}, projection-based nonintrusive MOR methods are interested in finding a subspace $\mathcal{V}_r$ and ROM
\begin{equation}\label{eq:SDEOIROM}
\rmd \hX_r(t) = \left[\hA_r\hX_r(t) + \hB_ru(t) + \sum_{i=1}^m \hN_{r,i}\hX_r(t)u_i(t)\right]\rmd t + \hM_r \rmd \hat W(t),
\end{equation}
with $\hX_r(0) = X_{r,0}$, the coefficients $\hA_r,\hN_{r,1},\dots,\hN_{r,d}\in\R^{r\times r}$, $\hB_r\in\R^{r\times m}$, $\hM_r\in\R^{r\times d_r}$ and a $d_r$-dimensional Wiener process $\hat W$, such that $\V_r\hX_r(t)\approx X(t)$ as well. However, this is achieved without explicit access to the FOM drift coefficients $\A,\B, \N_1,\dots,\N_m$, the diffusion coefficient $\M$ or the correlation matrix $\K$ of the noise generating process. As $\hat W$ generally differs from $W$, it is important to notice that the approximation by \eqref{eq:SDEOIROM} cannot be meant in a pathwise or $L^2$-sense which is in contrast to previous MOR approaches for SDEs.
The recently developed operator inference (OpInf) method \cite{peherstorfer_data-driven_2016,Ghattas_Willcox_2021,Kramer2024Review} provides an approach for ODEs with bilinear terms and nonlinearities of polynomial structure by solving a least squares minimisation problem that fits reduced coefficients to projected FOM state observations. 
Extensions to nonlinearities via variable transformations \cite{qian_lift_2020,qian_reduced_2022,khodabakhshi_non-intrusive_2022}, second order systems \cite{filanova_mechanical_2023}, as well as the combination with an interpolation method for non-polynomial nonlinearities with analytic form \cite{benner_operator_2020}, are available.
The case of noisy or low-quality data was considered by Uy et al. \cite{uy_operator_2023,uy_operator_2021,uy_active_2023}.
Guo et al. \cite{GUO2022115336} investigate the potential of a probabilistic interpretation of OpInf from the perspective of Bayesian inverse problems.
Studies on projection-based and interpolation methods can be found, for instance, in \cite{benner_survey_2015,antoulas_interpolatory_2020,benner_model_2017}.

The goal of this paper is to establish an OpInf approach for SDEs. This is an enormous challenge, since, in addition, the noise process cannot be observed, i.e., even if a path of the SDE solution is available, we do not know the associated trajectory of the driving process.
This also means that we do not know the dimension $d$ of the Wiener process.
Therefore, it is not possible to construct pathwise accurate approximations. 
Instead, our ansatz reproduces the distribution of the FOM state variable. 
We present the details of this work in several sections structured as follows. 
Section \ref{sec:MOR} introduces the grey box setting of bilinear SDEs under the influence of additive Gaussian noise. This means that the structure of the full order SDE is known, but not its coefficients. We only have access to full order state observations.
Building on the briefly explored distributional properties of the FOM, the inference methods of a reduced drift and diffusion coefficients are constructed.
A summary of the developed approach by Algorithm \ref{alg:OpInfSDE} concludes the section.
In Section \ref{sec:SubspaceConstruction}, the dominant subspace construction based on snapshots, known from POD, is generalised in two ways.
These generalisations arise from the way the snapshot matrix is constructed. 
To this end, the \emph{state snapshot matrix} and the \emph{moment snapshot matrix} are introduced.
The respective minimisation problems, concerning the optimal projection, are closely related by construction. 
Regarding the respective non-zero singular values, it is shown that, even for general SDEs, the span of the corresponding left singular vectors of the moment snapshot matrix is always contained in the span of the left singular vectors of the state-snapshot matrix.
Section \ref{sec:Closeness} establishes the closeness of the ROM obtained by the developed approach to the intrusive POD-ROM. 
To this end, Theorems \ref{thm:convergence_Expectation} and \ref{thm:convergence_Covariance} demonstrate that, if the training data is generated using consistent estimators for the mean and the covariance, the OpInf ROM coefficients converge almost surely to the POD-ROM coefficients. 
In the case where only the FOM state observations are available, Proposition \ref{pro:error_estimators} shows that the OpInf ROM coefficients can be arbitrarily close to the POD-ROM coefficients with probability $1$, if the number of training samples, the time step size and the ROM dimension are chosen appropriately.
Theorems \ref{thm:CovErrSol} and \ref{thm:ExpErrSol} present a bound for the mean and the covariance deviation of two ROMs depending on the distance of the respective ROM coefficients. 
Hence, under the conditions detailed in Proposition \ref{pro:error_estimators}, the mean and the covariance of the OpInf ROM and POD-ROM are arbitrarily close with probability $1$.
As the expectation and the covariance fully determine the distribution of the state in our setting, we can hence expect a very good approximation of the statistical properties of the FOM, provided enough data are collected.
Section \ref{sec:Experiments} provides two numerical experiments, which compare the developed approach to the POD method. 
This is done by using the summed relative errors in expectation and covariance, respectively, as well as the relative weak error with respect to two functionals, evaluated at the end-time $T$.

\section{Nonintrusive Model Order Reduction in the SDE setting}
\label{sec:MOR}
In the deterministic case of \eqref{eq:SDEFOM} ($\M\equiv 0$ or $\K\equiv 0$), the OpInf approach and the extensions thereof \cite{peherstorfer_data-driven_2016,benner_operator_2020, qian_lift_2020,uy_operator_2021,qian_reduced_2022,khodabakhshi_non-intrusive_2022,uy_active_2023,uy_operator_2023,filanova_mechanical_2023}, can provide an effective and easy to implement reduction method that solely requires the accessibility of observations of system-state, input $u$ and the initial condition.
The original method obtains the reduced coefficients by solving a least-squares problem, where the right-hand side is constructed from approximated time-derivatives of the state-trajectory at the observation times.

However, the stochastic case poses the challenge that the paths of $X(t),0\leq t\leq T$ are not differentiable, thereby prohibiting the application of this method. 
Furthermore, observations of the noise process $W$ are usually not available. 
Hence, a pathwise approximation, also known as strong approximation, is infeasible.
Therefore, the proposed method aims to obtain a ROM that approximates the FOM in the weak sense, i.e., that the approximation is with regard to the distribution of the FOM at each time $t\in[0,T]$. 

A general linear SDE \cite{kloeden_numerical_1992} can be solved by utilising the corresponding fundamental matrix. Considering the FOM \eqref{eq:SDEFOM}, one can find the explicit solution representation
\begin{equation}\label{eq:SDE_AnalyticalSolution}
X(t) = \Phi(t)\left( X_0 + \int_0^t \Phi(s)^{-1}\B u(s)\rmd s + \int_0^t \Phi(s)^{-1}\M\rmd W(s)\right),
\end{equation}
where $\Phi(t)\in\R^{n\times n}$ is the fundamental matrix that solves the differential equation
\begin{align*}
\dot{\Phi}(t) = \left[ \A + \sum_{i=1}^m \N_iu_i(t)\right] \Phi(t), \quad \Phi(0) = \I_n. 
\end{align*}
Here, $\I_n$ denotes the identity matrix with $n$ rows and columns. 
Since the control $u$ is chosen to be deterministic, the integral $\int_0^t \Phi(t)\Phi(s)^{-1}\M\rmd W(s)$ can be interpreted in the Wiener sense and is therefore a Gaussian random vector for each time $t\in[0,T]$.
Hence, $X(t)$ is also Gaussian and therefore determined, in distribution, by its expectation $E(t):=\E{X(t)}\in\Rn$ and covariance $C(t):=\E{(X(t)-E(t))(X(t)-E(t))^T}\in\Rnn$.
Furthermore, since the ROMs \eqref{eq:SDEROM} and \eqref{eq:SDEOIROM} are of the same structure as the SDE \eqref{eq:SDEFOM}, the distribution of the reduced variables is Gaussian as well. 
As mentioned above, observations of $W$ are generally unavailable. 
In the remainder of this section, an OpInf method is developed that solely requires samples of observations of the FOM system state, but does not rely on the observations of the paths of $W$.
As a consequence, the ROM is designed to approximate the FOM in distribution.
In the setting of \eqref{eq:SDEFOM} this equates to constructing a ROM with expectation $E_r(t)\in\R^r$ and covariance $C_r(t)\in\R^{r\times r}$, such that $\V_rE_r(t) \approx E(t)$ and $\V_rC_r(t)\V_r^T \approx C(t)$, rather than an approximation of the FOM system state directly. Here, $\V_r$ is a suitable matrix determining the projection.

\paragraph{Drift operator inference}
Using $E(t) = \E{X(t)}$, one obtains that the dynamics of the FOM expectation are governed by 
\begin{align}
E(t) &= E(0) + \int_0^t \A E(s) + \B u(s) + \sum_{i=1}^m \N_{i} E(s)u_i(s) \rmd s \nonumber\\
\iff
\dot{E}(t) &= \A E(t)+\B u(t)+ \sum_{i=1}^m \N_{i} E(t)u_i(t). \label{eq:MOR_ODE_ExpFOM}
\end{align} 
The statement is obtained by representing equation \eqref{eq:SDEFOM} in integral form, taking expectations on both sides and exploiting that the Itô integral has mean zero. 
Note that the noise generating Wiener process $W$ does not appear in the expectation dynamics.
The trajectory of $E(t), t\in[0,T]$, is therefore smooth almost everywhere if the control $u$ is integrable.
Similarly, one obtains an ODE of the same structure
\begin{equation}
    \dot{E}_r(t) = \A_r E_r(t)+\B_r u(t)+ \sum_{i=1}^m \N_{r,i} E_r(t)u_i(t) \label{eq:MOR_ODE_ExpROM}
\end{equation}
for the expected value $E_r(t):=\E{X_r(t)}$ of the POD-ROM \eqref{eq:SDEROM}.
Following the ideas of deterministic OpInf, let the $s+1$ observation times between $t_0=0$ and $t_s=T$ be denoted by $t_i,\, i=0,\dots,s$ and define 
\begin{subequations}
\label{eq:OIExpNotation}
\begin{align} 
\rmE_r &:= \begin{bmatrix} E_r(t_0),\dots,E_r(t_s)\end{bmatrix} \in{\R}^{r\times s+1},\\
\U &:= \begin{bmatrix} u(t_0),\dots, u(t_s) \end{bmatrix} \in{\R}^{m\times s+1}, \\
\U\odot\rmE_r &:=  \begin{bmatrix} u(t_0)\otimes E_r(t_0),\dots,u(t_s)\otimes E_r(t_s) \end{bmatrix} \in{\R}^{mr\times s+1},\\
\D &:= \begin{bmatrix} \rmE_r^T,\U^T,(\U\odot\rmE_r)^T \end{bmatrix}^T \in \mathbb{R}^{(r+m+mr) \times s+1},\\
\rmR &:= \begin{bmatrix} \dot{E}_r(t_0),\dots, \dot{E}_r(t_s)\end{bmatrix} \in{\R}^{r\times s+1}.
\end{align}
\end{subequations}
Exploiting the affine linear dependence of \eqref{eq:MOR_ODE_ExpROM} on $E_r(t)$ and $u(t)$, one can easily see that
\begin{align*}
\rmR = \O_r \cdot \D,
\end{align*}
 where $\O_r := \begin{bmatrix}
    \A_r,\B_r,\N_{r,1},\dots,\N_{r,m}
\end{bmatrix}\in\R^{r\times r+m+mr}$
is satisfied exactly.
The OpInf approach now consists of formulating and solving the inverse problem 
\begin{equation}\label{eq:Drift_LSQ}
\hO_r = \argmin_{\tO_r\in\R^{r\times (r+m+mr)}} \|\rmR^T-\D^T\tO_r^T\|_F^2
\end{equation}
for given $\D$ and $\rmR$, e.g., based on observations.
Furthermore, if the data-matrix $\D$ has full row-rank, the solution is unique.

However, direct observations of $\dot{E}_r(t)$ or of $E_r(t)$ are generally unavailable.
Instead, only finitely many samples of the state at discrete times $t_0=0,\dots,t_s=T$ are available.
The \emph{state snapshot matrix} is defined to be the column-wise arrangement of all $L$ sampled realisations of the states of \eqref{eq:SDEFOM} at the time steps $t_0,\dots,t_s$ 
\begin{equation}\label{eq:Snapshots}
\bbX := \begin{bmatrix}
    X(t_0,\omega_1),\dots,X(t_0,\omega_L,)\dots,X(t_s,\omega_1),\dots, X(t_s,\omega_L)
\end{bmatrix}\in\R^{n\times (s+1)L}.
\end{equation}
Note that the observations $X(t_i,\omega_j)$ might be obtained by numerical simulations and therefore are samples of the state of time-discretised dynamical system, rather than samples of the time-continuous solution of the FOM SDE. 
Hence, any approximation of the moments of $X(t)$ by such observations is influenced by two error sources: the finite number of samples $L$ and the time-discretisation of the continuous dynamics.
The columns $\rmE_{r,i}^L,i=0,\dots,s$ of the empirical mean of the projected states
\begin{equation}\label{eq:DataMean}
\rmE_{r}^L :=
\begin{bmatrix}
\rmE_{r,0}^L,\dots,\rmE_{r,s}^L
\end{bmatrix}, \text{ with }
\rmE_{r,i}^L = \frac{1}{L} \sum_{j=1}^L \V_r^TX(t_i,\omega_j)
\end{equation}
approximate the trajectory of the expected value of the reduced trajectory $\E{\V_r^TX(t_i)}$ at the observation times $t_0=0,\dots,t_s=T$. 
Generally, computing our reduced system \eqref{eq:SDEOIROM} is not equivalent to projecting the FOM resulting, e.g., in the POD model. The POD approach collects data from the Markov-process $X$ in order to find $\V_r$ leading to a reduced order Markov-process $X_r$ solving \eqref{eq:SDEROM}. The solution $\hX_r$ of the data-based model \eqref{eq:SDEOIROM} aims to mimic $X_r$. However, we do not have access to the Markov-process $X_r$ for estimating the coefficients of the system. Instead, we solely have data from $\V_r^T X$, which generally is not Markovian.
Hence, since $\V_r^TX(t_i)\neq X_r(t_i)$ in general and consequently $E_r(t_i)\neq \E{\V_r^TX(t_i)}$, as well as $C_r(t_i)\neq\Cov(\V_r^TX(t_i),\V_r^TX(t_i))$ in the considerations below,
this introduces an additional error into the training data, which does not vanish with increasing sample size.
Here, we have denoted the covariance matrix of two random vectors $X,Y$ of the same dimension with $\Cov(X,Y):=\E{\left( X-\E{X}\right)^T\left( Y-\E{Y}\right)}$.
This error can be reduced by increasing the ROM dimension $r$ or, under certain conditions, be eliminated by adjusting the data sampling scheme \cite{Peherstorfer2020Preasymtotic,uy_operator_2021}.
The precise choice of the projection matrix $\V_r$ depends on the particular approach that we specify in Section \ref{sec:SubspaceConstruction}.
If the time between the observations is sufficiently small, the time derivative of $\dot{E}_r(t_i), i=0,\dots,s$ can be approximated from $\rmE_{r}^L$, for example, by a finite difference scheme.
The result of such an approximation of $\dot{E}_r(t)$ at the time $t=t_i$ is denoted with $\dot{\rmE}^{L,h}_{r,i}$.
Now, utilising $\rmE_{r}^L$ and $\dot{\rmE}^{L,h}_{r}:=[\dot{\rmE}^{L,h}_{r,0},\dots, \dot{\rmE}^{L,h}_{r,s}]$, instead of $E_{r}(t)$ and $\dot{E}_{r}(t)$, to construct 
\begin{align*}
    \D^{L,h} &= \begin{bmatrix} (\rmE_{r}^L)^T,\U^T,(\U\odot\rmE_{r}^L)^T \end{bmatrix}^T \in \mathbb{R}^{(r+m+mr) \times s+1} \text{ and }\\
    \rmR^{L,h} &:= \begin{bmatrix} \dot{\rmE}^{L,h}_{r,0},\dots, \dot{\rmE}^{L,h}_{r,s}\end{bmatrix} \in{\R}^{r\times s+1}.
\end{align*}
One can solve the perturbed version of the least squares problem \eqref{eq:Drift_LSQ}
\begin{equation}
    \label{eq:Drift_LSQ_practical}
   \hO^{L,h}_r = \argmin_{\tO_r\in\R^{r\times (r+m+mr)}} \|(\rmR^{L,h})^T-(\D^{L,h})^T\tO_r^T\|_F^2 
\end{equation}
to obtain an approximation $\hO^{L,h}_r=[\hA^{L,h}_r,\hB^{L,h}_r,\hN^{L,h}_{r,1},\dots,\hN^{L,h}_{r,m}]$ of $\hO_r$.
If one can guarantee that $\D^{L,h}$ and $\D$ have full rank, then $\hO^{L,h}_r$ is an estimate of $\O_r$, since the solutions of \eqref{eq:Drift_LSQ} and \eqref{eq:Drift_LSQ_practical} are unique.

\paragraph{Diffusion operator inference}
Similar as before, the time evolution of the covariance matrix $C(t)$ is described by a (Lyapunov) differential equation.
First, the centred process $X^C(t) := X(t) - E(t)$ satisfies the SDE
\[\rmd X^C(t) = \left[\A X^C(t)+\sum_{i=1}^m \N_{ i} X^C (t)u_i(t)\right]\rmd t + \M \rmd W(t).\]
Next, applying Itô's product rule to the product $X^C (t)X^{C} (t)^T$ results in
\begin{align*}
\rmd (X^C (t)X^{C } (t)^T)
&=[ \A X^C (t)X^{C } (t)^T+X^C (t)X^{C } (t)^T \A^T  \\
&\quad + \sum_{i=1}^m\N_{ i}u_i(t)X^C (t)X^{C } (t)^T  + \sum_{i=1}^m X^C (t)X^{C } (t)^T (\N_{ i} u_i(t))^T ] \rmd t\\
&\quad + \M \rmd W(t)X^{C } (t)^T + X^{C} (t)\rmd W (t)^T\M^T  + \M \K\M^T \rmd t.
\end{align*}
Again, due to inaccessibility of the noise realisations, expectations are taken on both sides. 
With $C (t)=\E{X^{C} (t)X^{C } (t)^T}$, the zero mean of the Itô integral and slight algebraic manipulations, one obtains
\begin{equation*}
\frac{\rmd}{\rmd t} C (t) = \Psi(t) C(t)  + C (t) \Psi(t)^T+ \M \K\M^T,
\end{equation*}
where $\Psi(t) := \A + \sum_{i=1}^m \N_{ i}u_i(t)$ and the initial condition is given by $C(0) =\Cov(X_0,X_0)$. 
For the POD-ROM one obtains a similar system
\begin{equation*}
\frac{\rmd}{\rmd t} C_r (t) = \Psi_r(t) C_r(t)  + C_r (t) \Psi_r(t)^T+ \M_r \K\M_r^T,
\end{equation*}
with $\Psi_r(t) := \A_r + \sum_{i=1}^m \N_{r, i}u_i(t)$ and $C_r(0) =\Cov(X_{r,0},X_{r,0})= \V_r^T \Cov(X_0,X_0)\V_r$.
The influence of the inhomogeneities $\H =\M\K\M^T$ and $\H_r =\M_r\K\M_r^T$ on the previous equations is independent of time and the corresponding drift coefficients.
Defining $\S_{r,i}$ as the residual 
\begin{equation}\label{eq:CovResidual}
\S_{r,i} := \dot{C}_r(t_i)- \left[\Psi_r(t_i)C_r(t_i) + C_r(t_i)\Psi_r(t_i)^T\right],
\end{equation}
the least squares problem 
\begin{equation}\label{eq:MOR_CovMinimisation}
\hH_r = \argmin_{\tilde{\H}\in \R^{r\times r}} \sum_{i=1}^{s} \|\S_{r,i} - \tilde{\H}\|_F^2
\end{equation}
is solved by $\H_r$.
To obtain a practically feasible minimisation problem, the (reduced) empirical covariances 
\begin{align}\label{eq:DataCovariance}
\C_{r,i}^L:= \frac{1}{L-1}\sum_{k=1}^L (\V_r^TX(t_i,\omega_k)-\rmE_{r,i}^L)(\V_r^TX(t_i,\omega_k)-\rmE_{r,i}^L)^T,\quad i=0,\dots,s
\end{align}
are computed from the collected snapshot data $\bbX$ and involve a projection matrix $\V_r$ specified in Section \ref{sec:SubspaceConstruction}.
The approximations $\dot{\C}_{r,i}^{L,h}$ of the time derivatives are obtained, for instance, again by a finite difference scheme. 
Subsequently, a perturbed version 
\begin{equation}\label{eq:Cov_LSQ_practical}
\hH_r^{L,h} = \argmin_{\tilde{\H}\in \R^{r\times r}} \sum_{i=1}^{s} \|\hS_{r,i}^{L,h} - \tilde{\H}\|_F^2,
\end{equation}
of \eqref{eq:MOR_CovMinimisation} can be constructed using
\begin{align*}
    \hPsi_{r,j}^{L,h} &:= \hA_r^{L,h} + \sum_{i=1}^m \hN_{r, i}^{L,h}u_i(t_j),\\
    \hS_{r,i}^{L,h} &:= \dot{\C}_{r,i}^{L,h} - \left[\hPsi_{r,i}^{L,h}\C_{r,i}^L + \C_{r,i}^L(\hPsi_{r,i}^{L,h})^T\right],
\end{align*}
solely from the available snapshot data and already inferred drift coefficients.
Finally, one needs to obtain the diffusion coefficient and correlation matrix of the noise process of the surrogate model from $\hH_r^{L,h}\in\R^{r\times r}$.
This can be achieved by decomposing $\hH_r^{L,h}$ into a product $\hM_r^{L,h}\hK_r^{L,h}(\hM_r^{L,h})^T$, with $\hM_r^{L,h}\in\R^{r\times d_r}$ and $\hK_r^{L,h}\in\R^{d_r\times d_r}$ symmetric positive definite.
As the dimension $d$ of the FOM noise process is unknown, the dimension of the surrogate noise process $d_r$ has to be inferred from the data, for instance, by using the number of positive eigenvalues of $\hH_r^{L,h}$.
The noise generating process of the surrogate model is then constructed as a Wiener process with correlation matrix $\hK_r^{L,h}$.
Algorithm \ref{alg:OpInfSDE} summarises the approach detailed above.
\begin{algorithm}
\caption{operator inference for SDE \eqref{eq:SDEFOM}}
\label{alg:OpInfSDE}
\begin{algorithmic}[1]
\STATE{Collect observations $X(t_0,\omega_j),\dots, X(t_s,\omega_j), j=1,\dots,L$.}
\STATE{Choose projection matrix} $\V_r\in\R^{n\times r}$.
\STATE{Construct $\D^{L,h}$ and $\rmR^{L,h}$.}
\STATE{Solve the least squares problem \eqref{eq:Drift_LSQ_practical}} for $\hO^{L,h}_r$.
\STATE{Extract system coefficients $\hA_r^{L,h},\hB_r^{L,h},\hN_{r,1}^{L,h},\dots,\hN_{r,m}^{L,h}$ from $\hO^{L,h}_r$.}
\STATE{Construct $\hPsi^{L,h}_{r,i}$ and $\hS_{r,i}^{L,h}$.}
\STATE{Solve least squares problem \eqref{eq:Cov_LSQ_practical}} for $\hH_r^{L,h}$.
\STATE{Factorise $\hH_r^{L,h}$ into $\hM_r^{L,h}\in\R^{r\times d_r}$ and $\hK_r^{L,h}\in\R^{d_r\times d_r}$. }
\STATE Return $\hA_r^{L,h},\hB_r^{L,h},\hN_{r,1}^{L,h},\dots,\hN_{r,m}^{L,h},\hM_r^{L,h},\hK_r^{L,h}$.
\end{algorithmic}
\end{algorithm}

\paragraph{Practical considerations}
The inference from finite data poses several practical problems.
For instance, although $\hH_r$ is symmetric and positive (semi)definite, this is generally not the case for $\hH_r^{L,h}$.
Hence, factorisation algorithms such as a Cholesky decomposition, which rely on these properties, cannot be directly applied to compute the decomposition $\hM_r^{L,h}\hK_r^{L,h}(\hM_r^{L,h})^T = \hH_r^{L,h}$.
To address this issue, the matrix $\hH_r^{L,h}$ is symmeterised before an eigenvalue decomposition
\[\frac{\hH_r^{L,h}+(\hH_r^{L,h})^T}{2} = \U_H\S_H\U_H^T\]
is computed.
By construction, the diagonal entries of $\S_H$ approximate the eigenvalues of $\M_r\K\M_r^T$.
Let $\lambda_{H,i},\, i=1,\dots,r$ denote the diagonal entries of $\S_H$.
To remove the components that correspond to noise, the eigenvalues below a threshold of 0.1\% of the maximum eigenvalue are truncated, leaving $d_r$ significant eigenvalues $\lambda_{i_1},\dots,\lambda_{i_{d_r}}$.
The surrogate diffusion coefficient is defined columnwise as $\hM_{r,j}^{L,h}:= \sqrt{\lambda_{i_j}}\U_{H,j}$ for all $j=1,\dots,d_r$.
The correlation matrix of the surrogate noise process is hence set to $\hK_{r}^{L,h} = \I_{d_r}$.
While alternative constructions are possible, this choice ensures that the surrogate model is driven by uncorrelated standard Gaussian noise, even if the FOM noise process $W$ has correlated components.

Imposing structural constraints as discussed above regularises the solution $\hH_r^{L,h}$ of \eqref{eq:Cov_LSQ_practical} by projecting onto a subset of all symmetric positive definite matrices.
Conversely, a standard technique to regularise the linear least-squares problem \eqref{eq:Drift_LSQ} is Tikhonov regularisation, which penalises the Frobenius norm of the inferred coefficients. 
The specific choice of regularisation depends heavily on the particular problem. 
If, for instance, the linear and bilinear coefficients $\A$ and $\N = [\N_{1},\dots,\N_{m}]$ have a different scale, one can solve
\[\hO_r = \argmin_{\substack{\tO_r\in\R^{r\times (r+m+mr)}\\ \tO_r = [\tilde \A_r,\tilde \B_r, \tilde \N_r]}} \|\rmR^T-\D^T\tO_r^T\|_F^2 + \gamma_1\left(\Vert \tilde\A_r\Vert_F^2 + \Vert \tilde \B_r\Vert_F^2\right) + \gamma_2\Vert \tilde\N_r\Vert_F^2\]
with appropriate $\gamma_1,\gamma_2>0$ instead of solving \eqref{eq:Drift_LSQ} to ensure that this scaling difference is present in the inferred coefficients as well.
Since the scaling parameters balance the fit of the regression and the imposed penalisation term, choosing a suitable set of regularisation parameters is not trivial. 
An adaptive approach of finding $\gamma_1$ and $\gamma_2$ is proposed in \cite{McQuarrie03042021} for a deterministic system, where a quadratic term replaces the bilinearity.
The method consists of minimising the ROM error over all regularisation parameters using a derivative-free method.
As OpInf is inherently a data-driven method, overfitting poses a problem, which can be mitigated by regularisation. 
Regularisation in the context of OpInf has been explored with great success \cite{peherstorfer_data-driven_2016,qian_lift_2020,qian_reduced_2022,McQuarrie03042021,Swischuk2020Learning} to ensure ROM stability and well-posedness in setting with limited data.

\section{Construction of $\mathcal{V}_r$}
\label{sec:SubspaceConstruction}
This section is concerned with the construction of an orthonormal matrix 
$\V_r\in\R^{n\times r}$ whose column vectors form a basis of $\mathcal{V}_r$. A projection $\rmP$ onto $\mathcal{V}_r$ can consequently be defined by $\rmP=\V_r\V_r^T\in\Rnn$.
In the classical deterministic POD setting, one chooses the subspace $\mathcal{V}_r$ and the corresponding projection $\rmP$ in such a way that $\rmP \bbX$ is the best low-rank approximation of rank $r$ of the snapshots $\bbX$ over all possible subspaces of dimension $r$.
That is, one computes the leading $r$ left singular vectors of $\bbX$ onto which the data is subsequently projected. 
The matrix $\V_r$ hence consists of the leading $r$ left singular vectors of $\bbX$ or, equivalently, the leading $r$ eigenvectors of $\bbX\bbX^T$.

However, to achieve a good approximation in the weak sense of our setting, ideally one could identify a subspace $\mathcal{V}_r\subset\Rn$ that is simultaneously optimal for the projection of the covariance and expectation, rather than the states of the FOM. 
To this end, two methods of constructing an appropriate snapshot matrix are proposed. 
Both methods reduce to the same POD snapshot matrix in the absence of noise.
If one collects the $L$ independent samples, the immediate generalisation of the POD method suggests arranging the observations into a \emph{state snapshot matrix} $\bbX$ as defined in \eqref{eq:Snapshots}, or a version thereof with permuted columns.
The projection defined by the leading $r$ left singular vectors of this state snapshot matrix corresponds to finding the dominant features of the averaged non-centralised second moments of $X(t,\omega),\,t\in\{t_0,\dots,t_s\}$.
This can be seen by observing that 
\begin{align*}
    \frac{1}{(s+1)(L-1)} \bbX\bbX^T &= \frac{1}{(s+1)(L-1)}\sum_{i=0}^{s}\sum_{j=1}^L X(t_i,\omega_j)X(t_i,\omega_j)^T\\
    &\rightarrow \frac{1}{s+1}\sum_{i=0}^s C(t_i) +  E(t_i)E(t_i)^T,
\end{align*}
for $L\to \infty$ by the law of large numbers. 
Using the non-centralised quadratic forms $\hat{\Z}:= \frac{1}{s+1}\sum_{i=0}^s C(t_i)+ E(t_i)E(t_i)^T$, averaged over the observation times, one obtains
\begin{align*}
   \frac{1}{(s+1)^2(L-1)^2} \| \bbX\bbX^T - \rmP\bbX\bbX^T\rmP \|_F^2 
    &\rightarrow \| \hat{\Z} - \rmP \hat{\Z} \rmP \|_F^2 
\end{align*}
for $L\to\infty$.
Since $\hat{\Z}\in\Rnn$ is symmetric and positive (semi-)definite,
one can compute an eigenvalue decomposition of the form
$\hat{\Z}= \U\Sigma\U^T,$ with $ \U \in\Rnn$ and $ \Sigma = \diag(\lambda_1,\dots,\lambda_n)\in\Rnn$ where the eigenvalues are ordered according to their magnitude $\lambda_1\geq\dots\lambda_n\geq 0$.
It is well known that the solution of minimisation of $\| \hat{\Z} - \rmP \hat{\Z} \rmP \|_F^2 $ over the possible orthogonal projections $\rmP$ of rank $r$ can be obtained by defining $\rmP$ as the projection onto the subspace spanned by the $r$ leading eigenvectors. 
In particular, this means that we define $\mathcal{V}_r=\mathrm{span}\{v_1,\dots,v_r\}$, $\V_r:=[v_1,\dots,v_r]$, and $\rmP := \V_r\V_r^T$, where $v_r$ is the eigenvector corresponding to the $r$-th largest eigenvalue $\lambda_r$.
The error of this projection then is $\| \hat{\Z} - \rmP \hat{\Z} \rmP \|_F^2=\textstyle\sum_{i=r+1}^n \lambda_i^2$.

However, since we are interested in the approximation in the weak sense, i.e., in distribution, one could consider a subspace $\mathcal{V}_r$ that best describes not the states of the FOM system, but the states of expectation $E(t)$ and covariance $C(t)$.
To obtain an orthogonal basis for the image space of both functions, one can investigate the column space spanned by the \emph{moment snapshot matrix}
\[\F:= [E(t_0),\dots, E(t_s),C(t_0),\dots,C(t_s)] 
= [\rmE,\C] \in\R^{n\times (n+1)(s+1)}\]
with $\rmE = [E(t_0),\dots, E(t_s)]\in\R^{n\times s+1}$ and 
$\C= [C(t_0),\dots,C(t_s)]\in\R^{n\times n(s+1)}$.
One can then compute an approximation of the orthogonal basis of the column space of $\F$ simply by computing the left singular vectors of the empirical moment snapshot matrix 
\[\F^L:= [\rmE_{f,0}^L,\dots, \rmE_{f,s}^L,\C_{f,0}^L,\dots,\C_{f,s}^L] 
= [\rmE^L,\C ^L] \in\R^{n\times (n+1)(s+1)},\]
where $\rmE^L_{f,i}$ and $\C^L_{f,i}$ are the empirical expectation and covariance of the full-state observations $X(t_i,\omega_j),$ $j=1,\dots,L$ at the time $t_i$.
The matrices $\rmE^L\in\R^{n\times s+1}$ and $\C^L\in\R^{n\times n(s+1)}$ are the respective matrices constructed from the column-wise stacking of the estimated moments.
The projection defined by the $r$ leading left singular vectors of $\F^L$ then solves the task of minimising
\[\Vert \F^L(\F^L)^T-\rmP\F^L(\F^L)^T\rmP \Vert_F^2 = \Vert \tilde\Z-\rmP\tilde\Z\rmP \Vert_F^2,\]
with $\tilde\Z = \sum_{i=0}^s \rmE_{f,i}^L(\rmE_{f,i}^L)^T+(\C^L_{f,i})^2$ over all possible orthogonal projections of rank $r$.
Note, that the quadratic form of the moment snapshot matrix is closely related to that of the state snapshots, with the difference being the omitted scaling by $(s+1)^{-1}$ and the quadratic appearance of the covariance. 
The remainder of this section briefly investigates the connection between the left singular vectors of the snapshot matrices.

Even outside the presented linear SDE with additive Gaussian noise setting, it can be established that subspace constructions by the left singular vectors of the state snapshots and moment-snapshots are closely related.
For the remainder of this section, let the FOM stochastic process $X(t,\omega)\in\Rn$ be defined by the SDE
\[\rmd X(t) = \mu(t,X(t))\rmd t + \sigma(t,X(t))\rmd W(t),\quad t\in [0,T].\]
The coefficient functions $\mu:[0, T]\times \Rn\to \Rn$ and $\sigma:[0, T]\times\Rn\to\R^{n\times d}$ are nice enough, so that the existence and uniqueness of a global solution $X$ is guaranteed. Hence, the expectation, and the covariance, defined by $ E(t):=\E{X(t)}$ and $C(t):=\E{(X(t)-E(t))(X(t)-E(t))^T}$, exist for all $t\in[0,T]$.
The restriction to $t\in[0,T]$ is not necessary and can be replaced by 
$t\in[0,\infty )$ without a change in the arguments.
The proofs in the remainder of this section make extensive use of the compact SVD (cSVD). 
That is, any decomposition of a rank $r_a$ matrix $\A = \V_{r_a}\S_{r_a}\W_{r_a}^T \in\Rnm$, such that $\V_{r_a}\in \R^{n\times r_a}$ and $\W_{r_a}\in\R^{m\times r_a}$ have orthonormal columns and $S_{r_a}\in\R^{r_a\times r_a}$ is a diagonal matrix with the non-zero singular values $\sigma_1\geq\dots\geq\sigma_{r_a}>0$ on the main-diagonal.
The following Theorem shows that the space spanned by the left singular vectors of $\rmE^L$ corresponding to non-zero singular values is a subspace of the span of the left singular vectors of $\bbX$.
\begin{theorem}\label{thm:LSV}
Let $L,s\geq 1$ and let the state-snapshot and empirical expectation snapshot matrices have the cSVDs 
\begin{align*}
 \bbX &= \V_r\S_r\W_r^T \text{ and } \\
 \rmE^L &=\V_{e,r_e}\S_{e,r_e}\W_{e,r_e}^T,
\end{align*}
respectively.
There exists a matrix $\T_e\in\R^{r\times r_e}$ such that $\V_{e,r_e}=\V_r\T_e$.
\end{theorem}
\begin{proof}
Without loss of generality, let
\[\bbX = \begin{bmatrix}
    X(t_0,\omega_1),\dots,X(t_0,\omega_L,)\dots,X(t_s,\omega_1),\dots, X(t_s,\omega_L)
    \end{bmatrix} \in\R^{n\times L(s+1)}\]
 be the ordering of the columns of the state snapshot matrix.
If $\rmP\in\R^{L(s+1)\times (s+1)}$ is defined by
\begin{align*}
    \rmP := \frac{1}{L} \I_{s+1}\otimes\mathbbm{1}_L, 
\end{align*}
with $\mathbbm{1}_L=[1,\dots,1]^T\in\R^L$,
then $\rmE^L=\bbX\rmP$, since
\begin{align*}
    \bbX\rmP &= \frac{1}{L} 
    \begin{bmatrix}
        X(t_0,\omega_1),\dots, X(t_0,\omega_L),\dots, X(t_s,\omega_1),\dots, X(t_s,\omega_L)
    \end{bmatrix}
    \cdot (\I_{s+1}\otimes\mathbbm{1}_L)\\
    &= 
    \frac{1}{L}
    \begin{bmatrix}
    \sum_{j=1}^L X(t_0,\omega_j),\dots,\sum_{j=1}^L X(t_s,\omega_j)
    \end{bmatrix}\\
    &=
    \begin{bmatrix}
    \rmE_{f,0}^L,\dots,\rmE_{f,s}^L
    \end{bmatrix}.
\end{align*}
Let $\bbX = \V_r\S_r\W_r^T$ and $\rmE^L = \V_{e,r_e}\S_{e,r_e}\W_{e,r_e}^T$ be the cSVDs of the state and expectation snapshot matrices $\bbX$ and $\rmE^L$, respectively.
Then,
\begin{align*}
    \rmE^L &= \bbX\rmP,\\
    \iff \V_{e,r_e}\S_{e,r_e}\W_{e,r_e}^T &= \V_r\S_r\W_r^T\rmP,\\
    \implies \V_{e,r_e} 
    &= \V_r\underbrace{\S_r\W^T_r\rmP\W_{e,r_e}\S^{-1}_{e,r_e}}_{=:\T_e},
\end{align*}
where $\S_{e,r_e}^{-1}=\mathrm{diag}(\sigma_{e,1}^{-1},\dots,\sigma_{e,r_e}^{-1})$ is the inverse of $\S_{e,r_e}$.
\end{proof}
Employing similar arguments as in the previous proof, analogous statements regarding the left singular vectors of the centralised snapshots and the covariance snapshots can be obtained. 
\begin{lemma}\label{lem:LSV3}
Let $L\geq 2,\, s\geq 1$ and let the matrix of empiric covariance snapshots have the cSVD $\C^L = \V_{c,r_c}\S_{c,r_c}\W_{c,r_c}^T$.
There exists a matrix $\T_c\in\R^{r\times r_c}$, such that $\V_{c,r_c} =\V_r\T_c$.
\end{lemma}

\begin{proof}
Notice that
    $
    \C^L = \begin{bmatrix}
        \C_{f,0}^L,\dots,\C_{f,s}^L
    \end{bmatrix} $
    can be written as 
    \begin{align*}
        \C^L = \frac{1}{L-1}\hat{\bbX}\diag(\hat\bbX_0^T, \dots,\hat\bbX_s^T),
    \end{align*}
    where $\hat{\bbX}$ is the matrix of centralised state snapshots
\[
    \hat{\bbX} = \begin{bmatrix}
        \hat{\bbX}_0,\dots,\hat{\bbX}_s
    \end{bmatrix} = \bbX - \rmE^L\otimes \mathbbm{1}_L^T
    \]
    with blocks 
    \[
    \hat{\bbX}_i = \begin{bmatrix}
        X(t_i,\omega_1)-\rmE_{f,i}^L,\dots,X(t_i,\omega_L)-\rmE_{f,i}^L
    \end{bmatrix}.
    \]
    Now, considering the cSVDs $\C^L = \V_{c,r_c}\S_{c,r_c}\W_{c,r_c}^T$, $\hat{\bbX} = \V_{z,r_z}\S_{z,r_z}\W_{z,r_z}^T$ and the previous result, one obtains
    \begin{align*}
        \V_{c,r_c}\S_{c,r_c}\W_{c,r_c}^T &= \C^L\\
        &= \frac{1}{L-1}\hat{\bbX}\diag(\hat\bbX_0^T, \dots,\hat\bbX_s^T)\\
        &= \frac{1}{L-1}\V_{r}\left(\S_r\W_r^T - \T_{e}\S_{e,r_e}\W_{e,r_e}^T\otimes \mathbbm{1}_L^T \right)\diag(\hat\bbX_0^T, \dots,\hat\bbX_s^T)\\
        \implies \V_{c,r_c} &= \V_r \T_c,        
    \end{align*}
where the transformation matrix $\T_c$ is
\begin{align*}
    \T_c & = \frac{1}{L-1} \left(\S_r\W_r^T - \T_{e}\S_{e,r_e}\W_{e,r_e}^T \otimes \mathbbm{1}_L^T\right)\diag(\hat\bbX_0^T, \dots,\hat\bbX_s^T)\W_{c,r_c}\S_{c,r_c}^{-1}.
\end{align*}
\end{proof}
Collecting the previous results, the following theorem shows that the subspace constructed from the left singular vectors of $\bbX$ is a superset of the subspace obtained from the empiric moment snapshot matrix $\F^L$.
\begin{theorem}
\label{thm:LSV4}
Let $L\geq2,\, s\geq 1$ and let $\V_{m,r_m}=[v_{m,1},\dots,v_{m,r_m}]\in\Rnn$ be the left singular vectors and $r_m$ the rank of the empiric moment snapshot matrix 
\[\F^L = \begin{bmatrix}
        \rmE^L,\C^L
    \end{bmatrix}
    =
    \begin{bmatrix}
        \rmE^L_{f,0},\dots,\rmE^L_{f,s},\C^L_{f,0},\dots,\C^L_{f,s}
    \end{bmatrix}.
\]
The space $\bar{\mathcal{V}}_{r_m} = \mathrm{span}\{v_{m,1},\dots,v_{m,r_m}\}$ is
a subspace of $\mathcal{V}_r=\mathrm{span}\{v_1,\dots,v_r\}$ spanned by the left singular vectors $\V_r=[v_1,\dots,v_r]$ of the state snapshot matrix $\bbX$ corresponding to the non-zero singular values.
\end{theorem}
\begin{proof}
Similarly to the previous statements, the formulation is equivalent to finding some matrix $\T_m\in\Rnn$, such that $\V_{m,r_m} = \V_r\T_m$.
Let $\V_{m,r_m}\S_{m,r_m}\W_{m,r_m}^T$ be a cSVD of the empiric moment snapshot matrix 
    $\F^L=\begin{bmatrix}
        \rmE^L,\C^L
    \end{bmatrix}$.
    The statement is proven by using the previous results
    \begin{align*}
        \V_{m,r_m}\S_{m,r_m}\W_{m,r_m}^T &= \begin{bmatrix}
            \rmE^L,\C^L
        \end{bmatrix}\\
        &= \begin{bmatrix}
            \V_{e,r_e}\S_{e,r_e}\W_{e,r_e}^T,\V_{c,r_c}\S_{c,r_c}\W_{c,r_c}^T
        \end{bmatrix}\\
        &=\V_r \begin{bmatrix}
            \T_e \S_{e,r_e}\W_{e,r_e}^T, \T_c \S_{c,r_c}\W_{c,r_c}^T
        \end{bmatrix}\\
        \implies \V_{m,r_m} &=\V_r\T_m,
    \end{align*}
    with 
    \begin{align*}
    \T_m &= \begin{bmatrix}
            \T_e \S_{e,r_e}\W_{e,r_e}^T, \T_c \S_{c,r_c}\W_{c,r_c}^T
        \end{bmatrix}\W_{m,r_m}\S_{m,r_m}^{-1}.
    \end{align*}
\end{proof}
Notice, that for a sufficiently large number of samples, the left singular vectors of $\F^L=\begin{bmatrix}    \rmE^L_{f,0},\dots,\rmE^L_{f,s},\C^L_{f,0},\dots,\C^L_{f,s}
\end{bmatrix}$ are approximately the left singular vectors of $\F=\begin{bmatrix}
    \rmE(t_0),\dots,\rmE(t_s),\C(t_0),\dots,\C(t_s)
\end{bmatrix}$.
The left singular vectors of $\bbX$ therefore provide an adequate basis for the dynamics of expectation and covariance. 
Hence, obtaining the projection basis from the SVD of $\bbX$ is not only numerically advantageous, but also faster and directly available, since the formation of $\rmE^L$ and $\C^L$ is omitted. 
One advantage of the subspace construction from the moment-snapshot matrix $\F^L$ instead of $\bbX$ is that it allows for a scaling of the considered moments, for instance, by normalisation, to balance the contributions of $\rmE^L$ and $\C^L.$
This might be beneficial, if the approximation of one of the moments is prioritized in the construction from $\bbX$. 
However, as $\bbX$ represents state observations, such a reweighting might result in a suboptimal subspace if one is interested in a strong approximation of the FOM, i.e., an approximation of the trajectories $X$. 
\section{Closeness to the POD-ROM}\label{sec:Closeness}
In this section, it is shown that under certain conditions, the ROM obtained by Algorithm \ref{alg:OpInfSDE} is close in distribution to the ROM obtained by POD. 
Theorems \ref{thm:convergence_Expectation} and \ref{thm:convergence_Covariance} show that, given consistent estimators for the mean and the covariance, the coefficients inferred from a finite number of samples converge almost surely to those obtained using direct observations of expectation and covariance.
Proposition \ref{pro:error_estimators} further explains that the estimators constructed from the projected state variable can approximate certain consistent estimators arbitrarily well, provided that the ROM dimension $r$ is chosen large enough.  Therefore, an appropriate selection of the number of samples $L$, the time step size $h$ and $r$ ensures that our approach provides coefficients close to the ones of the POD-ROM.
Finally, Theorems \ref{thm:CovErrSol} and \ref{thm:ExpErrSol} establish that the difference in expectation and covariance between two bilinear SDEs with additive noise is determined by the distance of their respective coefficients.

For this section, it is assumed that the choice of the subspace $\mathcal{V}_r$ is fixed. 
In the performed tests, a significant source of errors in the inference steps \eqref{eq:Drift_LSQ_practical} and \eqref{eq:Cov_LSQ_practical}, was the approximation of the time-derivative to obtain the respective right-hand sides $\rmR^{L,h}$ and $S_{r,i}^{L,h}$.
Since the empirical mean and covariance converge almost surely and in $L^2(\Omega, \mathcal F, \mathbb P)$ with a rate (arbitrarily) close to $\frac{1}{2}$, a large $L$ is required to obtain satisfactory approximations.
The influence of the noise is especially relevant in the perturbation of the right-hand side $\rmR^{L,h}$, where a Monte-Carlo and a time discretisation are combined. It turns out that $L$ has to be chosen even more carefully when $h$ is small. We illustrate this in the following remark.
\begin{remark}\label{rem:finite_differences_Lh}
Let us consider a scalar version of \eqref{eq:SDEFOM}. 
We approximate $\frac{\rmd}{\rmd t}\mathbb E[X(t)]$ by $\frac{1}{Lh}\sum_{i=1}^L X_i(t+h)-X_i(t)$, where $(X_i(t+h)-X_i(t))_{i=1, 2, \dots}$ are i.i.d. copies of $X(t+h)-X(t)$. The error is \begin{align*}
&\mathbb E \Big\vert \frac{\rmd}{\rmd t}\mathbb E[X(t)] - \frac{1}{Lh}\sum_{i=1}^L X_i(t+h)-X_i(t)\Big\vert^2\\
&= \Big\vert \frac{\rmd}{\rmd t}\mathbb E[X(t)] - \frac{\mathbb E[X(t+h)-X(t)]}{h}\Big\vert^2+\mathbb E \Big\vert\frac{\mathbb E[X(t+h)-X(t)]}{h}  - \frac{1}{Lh}\sum_{i=1}^L X_i(t+h)-X_i(t)\Big\vert^2\\
&= \Big\vert \frac{\rmd}{\rmd t}\mathbb E[X(t)] - \frac{\mathbb E[X(t+h)-X(t)]}{h}\Big\vert^2+\frac{1}{L}\Var \Big(\frac{X(t+h)-X(t)}{h}\Big),
\end{align*}
since the above mixed term vanishes in expectation. The first term is the error in the approximation of the derivative and can be made arbitrary small given $h$ is small enough. However, such a small $h$ might cause a large $\frac{X(t+h)-X(t)}{h}$, since $X$ is not differentiable, and hence a large variance in the second term. Think, e.g., of $X$ being a Wiener process yielding $\Var \Big(\frac{X(t+h)-X(t)}{h}\Big)=\frac{1}{h}$. This effect can only be compensated by a proper choice of $L$. Thus, decreasing the time between observations $h$, without increasing $L$, leads to a dominance of the noise.
\end{remark}
Remark \ref{rem:finite_differences_Lh} further shows that the error, which is introduced by using samples instead of direct observations, vanishes as the number of samples increases.
More generally, we can show that, in the presence of strongly consistent estimators, the inferred drift coefficients converge almost surely to those obtained by direct observations.
Moreover, under the conditions of the following theorem, the distance to the POD reduced coefficients can be made arbitrarily small.
\begin{theorem}\label{thm:convergence_Expectation}
Let $\bar E_{r,0}^L,\dots,\bar E_{r,s}^L$ be almost surely convergent estimators of the expected values $E_{r,0},\dots,E_{r,s}$. Suppose that $\dot{E}_{r,i}^h$ is the finite difference approximation of $\dot E_{r,i}$, where $h$ is the time step between the observation times. Let $\dot{\bar{E}}_{r,i}^{L,h}$ be the sampling based version of $\dot{E}_{r,i}^h$, where the exact expectations $ E_{r,0},\dots, E_{r,s}$ are replaced by $\bar E_{r,0}^L,\dots,\bar E_{r,s}^L$. 
\begin{enumerate}
    \item 
    Then, using $\bar\rmE^L_r := \begin{bmatrix}\bar E_{r,0}^L,\dots,\bar E_{r,s}^L, \end{bmatrix}$, we obtain 
        \begin{align*}
            \bar\D^L&:= \begin{bmatrix}(\bar\rmE_r^L)^T,\U^T, (\U\odot \bar\rmE_r^L)^T\end{bmatrix} \to \D = \begin{bmatrix}
                \rmE_r^T,\U^T, (\U\odot \rmE_r)^T
            \end{bmatrix}
        \end{align*}
        almost surely for $L\to \infty$.
    \item  Moreover, we have
    \begin{align*}
        \bar\rmR^{L,h} &:= \begin{bmatrix} \dot{\bar{E}}_{r,0}^{L,h},\dots, \dot{\bar{E}}_{r,s}^{L,h}\end{bmatrix} \to \rmR^h := \begin{bmatrix} \dot{E}_{r,0}^h,\dots, \dot{E}_{r,s}^h\end{bmatrix}
    \end{align*}
    almost surely for $L\to\infty$.
    \item If additionally $\D^T$ has full column rank, then the solution $\hat{\bar{\O}}^{L,h}_r$ of 
    \begin{align}\label{eq:LSQ_convergence}
   \hat{\bar{\O}}^{L,h}_r = \argmin_{\tO_r\in\R^{r\times (r+m+mr)}} \|(\bar\rmR^{L,h})^T-(\bar\D^{L})^T\tO_r^T\|_F^2 
    \end{align}
    can be arbitrarily close to the solution $\hO_r$ of the minimisation problem \eqref{eq:Drift_LSQ} with probability $1$, if one chooses large $L$ and a small enough $h$.
\end{enumerate} 
\end{theorem}

\begin{proof}
    The first statement follows by the construction of $\bar{\D}^L$.
    The second statement follows directly from the definition of the finite difference approximation scheme, for example,
    \begin{align*}
        \lim_{L\to\infty} \dot{\bar{E}}^{L,h}_{r,i} = \lim_{L\to\infty} \frac{1}{h}\left( \bar{E}^L_{r,i+1}-\bar{E}^L_{r,i}\right)
        = \frac{1}{h}\left( E_{r,i+1}-E_{r,i}\right)
        = \dot E_{r,i}^h.
    \end{align*}
    To prove statement three, we write $\bar\D^L$ and $\bar\rmR^{L,h}$ as perturbations of $\D$ and $\rmR$
    \begin{align*}
        \bar\D^L = \D + (\Delta \D)^L, \text{ and } \bar\rmR^{L,h} = \rmR + (\Delta \rmR)^{L,h},
    \end{align*}
    where $\Vert (\Delta \D)^L \Vert_2$ converges to 0 and $(\Delta \rmR)^{L,h}$ converges to the error of the finite difference approximation scheme $(\Delta \rmR)^{h}$ almost surely for $L\to\infty$.
    This implies that there exists a $L^*$ such that 
    $\Vert (\Delta D)^L \Vert_2$ is smaller than the smallest singular value of $\D$ for all $L\geq L^*$ with probability $1$. 
    Hence, the minimisation problem \eqref{eq:LSQ_convergence} can be solved exactly and the solution $\hat{\bar{\O}}^{L,h}_r$ is unique for all $L\geq L^*$.

    Column-wise perturbation analysis of the full rank least squares problem \cite[Theorem 5.3.1]{golub_matrix_2013} reveals that the relative error between $\hO_{r,i}$ and $\hat{\bar{\O}}_{r,i}^{L,h}$ can be bounded by 
    \begin{equation*}
    \| \hO_{r,i}-\hat{\bar{\O}}_{r,i}^{L,h}\|_2\leq C_i\delta_{L,h} + \mathcal{O}(\delta^2_{L,h}),
    \end{equation*}
    where $\delta_{L,h} = \max\left\{\frac{\| (\Delta \D)^{L}\|_2}{\| \D \|_2},\frac{\| (\Delta \rmR)^{L,h}\|_2}{\| \rmR \|_2}\right\}$ and some constant $C_i$.
    For any $\epsilon>0$ we can choose a sufficiently small $h$ and large enough $L$, such that $\delta_{L,h}<\epsilon$ with probability 1. 
    As this holds for all columns and $\epsilon>0$ is arbitrary, this implies that we can make the Frobenius norm of the difference between $\hat{\bar{\O}}_r^{L,h}$ and $\hO_r$ arbitrarily small by choosing appropriate $h$ and $L$.
\end{proof}
A similar result holds for the inference of the diffusion coeffcients as well.
\begin{theorem}\label{thm:convergence_Covariance}
    Let $\bar E_{r,0}^L,\dots,\bar E_{r,s}^L$ and $\bar C_{r,0}^L,\dots, \bar C_{r,s}^L$ be almost surely convergent estimators of the expected values $E_{r,0},\dots,E_{r,s}$ and the covariances $C_{r,0},\dots, C_{r,s}$, respectively.
    Suppose that $\dot C^h_{r,i}$ is the finite difference approximation of $\dot C_{r,i}$ and
    let $\db{C}^{L,h}_{r,i}$  be the sampling based version of $\dot C^h_{r,i}$, where the exact covariances are replaced by $\dot{\bar{C}}^L_{r,0},\dots,\dot{\bar{C}}^L_{r,s}$.
    Furthermore let the coefficient matricies $\bar\Psi^{L,h}_{r,j}$ and $\Psi^h_{r,j}$ be defined by
    \begin{align*}
        \bar\Psi^{L,h}_{r,j} :=\hb{\A}^{L,h}_r+\sum_{i=1}^m\hb{\N}^{L,h}_{r,i}u_i(t_j) \text{ and }
        \Psi^h_{r,j} :=\A^h_r+\sum_{i=1}^m\N^h_{r,i}u_i(t_j),
    \end{align*}
    where $\hb{\A}^{L,h}_{r},\hb{\N}^{L,h}_{r,i}$ are obtained from \eqref{eq:LSQ_convergence} and $\A^h_{r},\N^h_{r,i}$ are obtained from the version of \eqref{eq:Drift_LSQ}, where $\rmR$ is replaced by $\rmR^h$.
    Finally, assume that $\D^T$ has full column rank. 
    \begin{enumerate}
        \item Then, for the residuals
            \begin{align*}
                \bar\S_{r,i}^{L,h} :=& \dot{\bar{\C}}_{r,i}^{L,h} - \left[\bar\Psi_{r,i}^{L,h}\bar\C^{L,h}_{r,i} + \bar\C^{L,h}_{r,i}(\bar\Psi_{r,i}^{L,h})^T\right]\text{ and} \\
                \S^h_{r,i} :=& \dot{\C}^h_{r,i} - \left[\Psi^h_{r,i}\C^h_{r,i} + \C^h_{r,i}(\Psi^h_{r,i})^T\right] ,
           \end{align*}
           we have that $\bar\S_{r,i}^{L,h}\to \S^h_{r,i}$ almost surely for $L\to\infty$.
        \item The solution $\hb{\H}^{L,h}_r$ of
        \begin{align}\label{eq:LSQ_convergence_cov}
            \hat{\bar{\H}}_r^{L,h} = \argmin_{\tilde{\H}\in \R^{r\times r}} \sum_{i=1}^{s} \|\bar\S_{r,i}^{L,h} - \tilde{\H}\|_F^2
        \end{align}
    \end{enumerate}     
    can be arbitrarily close to the solution $\hH_r$ of the minimisation problem \eqref{eq:MOR_CovMinimisation} with probability $1$, if one chooses a large $L$ and small enough $h$.
\end{theorem}
\begin{proof}
    By Theorem \ref{thm:convergence_Expectation} we know that $\hb{\Psi}^{L,h}_{r,i}$ converges to $\Psi^h_{r,i}$ with probability $1$ for $L\to\infty$.
    Using the same arguments as in the proof of the previous Theorem, the sampling based finite difference approximation $\dot{\bar{\C}}^{L,h}_{r,i}$ converges almost surely to the finite difference approximation of $\dot\C_{r,i}$ from direct observations. 
    Hence, we can write the sample based quantities as 
    \[\bar\Psi^{L,h}_{r,i} = \Psi^h_{r,i} + (\Delta\Psi)^{L}_i \text{ and } 
    \db{\C}^{L,h}_{r,i} = \dot\C^h_{r,i} + (\Delta\dot\C)^{L}_i,
    \]
    where the norms of the perturbations converge to $0$ with probability $1$ as $L\to\infty$, which proves the first statement. 
    This means, for any $\epsilon>0$ we can find $L_{\epsilon}$ and $h_{\epsilon}$ such that for all $i=0,\dots,s$, the norm between the estimated residual $\bar\S^{L,h}_{r,i}$ and the exact residual
    \begin{align*}
     \S_{r,i} &= \dot{\C}_{r,i} - \left[\Psi_{r,i}\C_{r,i} + \C_{r,i}(\Psi_{r,i})^T\right]
    \end{align*}
    is smaller than $\epsilon$ with probability $1$.
    The second statement is proven by first rewriting minimisation problem \eqref{thm:convergence_Covariance} as
    \[
    \hb{\H}^{L,h}_r = \argmin_{\tilde{\H}\in\R^{r\times r}} \Vert\bar\S_{r}^{L,h}-(\mathbbm{1}_s\otimes\I)\tilde\H\Vert_F^2,
    \]
    using the properties of the Frobenius norm, where 
    \[\bar\S_{r}^{L,h} = \begin{bmatrix}
        (\bar\S_{r,0}^{L,h})^T,\dots,(\bar\S_{r,s}^{L,h})^T
    \end{bmatrix}^T,\]
    and observing that the norm of the perturbation of $\bar\S^{L,h}_r$ obeys the bound
    \begin{align*}
        \Vert \bar\S_{r}^{L,h} -\S_r^{L,h}\Vert_F^2 
        &= \Vert \begin{bmatrix} (\bar\S_{r,0}^{L,h})^T,\dots,(\bar\S_{r,s}^{L,h})^T \end{bmatrix}^T - \begin{bmatrix} \S_{r,0}^T,\dots,\S_{r,s}^T \end{bmatrix}^T\Vert_F^2\\
        &\leq (s+1)\max_{i=0,\dots,s}\Vert \bar\S^{L,h}_{r,i}-\S_{r,i}\Vert_F^2\\
        &\leq (s+1)\epsilon^2
    \end{align*}
with probability $1$.
The result follows by employing the same perturbation analysis as in Theorem \ref{thm:convergence_Expectation}.
\end{proof}
Theorems \ref{thm:convergence_Expectation} and \ref{thm:convergence_Covariance} tell us that, if strongly consistent estimators were available, the more samples we collect to estimate the mean as well as the covariance and the smaller we choose $h$ to approximate their derivatives, the closer we are to the ROM coefficients of the POD-ROM \eqref{eq:SDEROM}. Such estimators are, e.g., 
\begin{align}\label{example_estimaor}
 \bar \rmE_{r,i}^L = \frac{1}{L} \sum_{j=1}^L X_r(t_i,\omega_j),\quad
 \bar \C_{r,i}^L:= \frac{1}{L-1}\sum_{k=1}^L (X_r(t_i,\omega_k)-\bar \rmE_{r,i}^L)(X_r(t_i,\omega_k)-\bar\rmE_{r,i}^L)^T.   
\end{align}
In practise we are only able to observe the FOM state $X$ and hence can only approximate the $E_r$ and $C_r$ at time $t_i$ by using the estimators $\rmE^L_{r,i}$ and $\C^L_{r,i}$, introduced in \eqref{eq:DataMean} and \eqref{eq:DataCovariance}, respectively. However, as these estimators rely on the projected FOM state variable, they are not consistent.
As mentioned after \eqref{eq:DataMean}, this is because generally $\V_r^TX(t_i)\neq X_r(t_i)$, where $\V_r^TX$ is not even Markovian. This means that we have to control the errors $\rmE^L_{r,i}-\bar \rmE^L_{r,i}$ and $\C^L_{r,i}-\bar \C^L_{r,i}$. If we can make it sufficiently small, then Algorithm \ref{alg:OpInfSDE} can provide coefficients close to the one of the POD-ROM. It turns out that a proper choice of $L$ and $h$ is not enough. We can only ensure closeness if the reduced dimension $r$ is selected appropriately. Below, we state a proposition pointing out this aspect.
\begin{proposition}\label{pro:error_estimators}
Let the ROM dimension $r$ so that $\Vert X(t)-\V_rX_r(t)\Vert<\delta$ almost surely for all $t\in [0, T]$ and a given $\delta>0$. Suppose that $\rmE^L_{r,i}$ and $\C^L_{r,i}$ are like in \eqref{eq:DataMean} and \eqref{eq:DataCovariance}. If $\bar \rmE^L_{r,i}$ and $\bar \C^L_{r,i}$ are given like in \eqref{example_estimaor}, then there exist two constants $c_1,c_2>0$, such that 
\begin{align*}
 \Vert \rmE^L_{r,i}-\bar \rmE^L_{r,i}\Vert<\delta \quad \text{and}\quad \Vert\C^L_{r,i}-\bar \C^L_{r,i}\Vert < c_1\delta + c_2\delta^2
\end{align*}
for all $i=0, \dots, s$ almost surely.
\end{proposition}
\begin{proof}
 Using that $\V_r$ is a matrix with orthonormal columns, we have \begin{align*}
     \Vert \rmE^L_{r,i}-\bar \rmE^L_{r,i}\Vert &= \frac{1}{L} \Vert\sum_{j=1}^L  (\V_r^T X(t_i,\omega_j)- X_r(t_i,\omega_j) )\Vert \\
     &\leq \frac{1}{L} \sum_{j=1}^L  \Vert \V_r^T X(t_i,\omega_j)- \V_r^T \V_r X_r(t_i,\omega_j) \Vert \\
     &\leq \frac{1}{L} \sum_{j=1}^L  \Vert X(t_i,\omega_j)- \V_r X_r(t_i,\omega_j) \Vert <\delta
 \end{align*}   
applying our assumption.
Computing the difference between $\bar\C^L_{r,i}$ and $\C^L_{r,i}$ is more involved.
To prove the statement, we employ the respective quantities again as perturbations
\begin{align*}
    \V_r^TX(t_i,\omega_j) &= X_r(t_i,\omega_j) + (\Delta X)_{i,j} \text{ and } \\
    \rmE^L_{r,i} &= \bar\rmE^L_{r,i} + (\Delta E)_i,
\end{align*}
where $\Vert (\Delta X)_{i,j} \Vert < \delta $ and $\Vert(\Delta E)_{i}\Vert<\delta$ almost surely.
We can express the empiric covariance matrices as
\begin{align*}
    \C_{r,i}^L &= \frac{1}{L-1} \sum_{j=1}^L \V_r^TX(t_i,\omega_j) X(t_i,\omega_j)^T\V_r + \frac{L}{L-1} \rmE^L_{r,i}\rmE^{L,T}_{r,i},\\
    \bar\C_{r,i}^L &= \frac{1}{L-1} \sum_{j=1}^L X_r(t_i,\omega_j) X_r(t_i,\omega_j)^T + \frac{L}{L-1} \bar\rmE^L_{r,i}\bar\rmE^{L,T}_{r,i}.
\end{align*}
Bounding the norm of the difference between $\C^L_{r,i}$ and $\bar\C^L_{r,i}$ can hence be reduced to bounding the norms of the differences between the respective summands.
Note, that for any four matrices $\A,\B,\C,\D$ of suitable dimensions, we have that 
\[\A\B^T - \C\D^T = \left(\A-\C\right)\left(\B-\D\right)^T + \C\left(\B-\D\right)^T + \left(\A-\C\right)\D^T. \]
First, using $\A=\B = X_r(t_i,\omega_j)$ and $\C=\D=\V_r^TX(t_i,\omega_j)$, we obtain
\begin{align*}
    \Vert \V_r^TX(t_i,\omega_j)X(t_i,\omega_j)^T\V_r - X_r(t_i,\omega_j)X_r(t_i,\omega_j)^T\Vert \leq \delta^2 + 2\Vert X(t_i,\omega_j)\Vert\delta ,
\end{align*}
since $\Vert\V_r^Tx\Vert\leq\Vert \V_r^T\Vert\Vert x\Vert=\Vert x\Vert$ for any $x\in\Rn$.
To express the bound on the difference between the empirical expectation estimators we can use
\begin{alignat*}{2}
    \A &= X_r(t_i,\omega_{j_1}), \quad \quad \B &&= X_r(t_i,\omega_{j_2}) \text{ and }\\
    \C &= \V_r^TX(t_i,\omega_{j_1}), \quad\D &&= \V_r^TX(t_i,\omega_{j_2}),
\end{alignat*}
where $j_1,j_2\in\{1,\dots,L\}.$
Expressing the empiric estimators by their respective definitions results in
\begin{align*}
    \Vert \rmE^L_{r,i}\rmE^{L,T}_{r,i} - \bar\rmE^L_{r,i}\bar\rmE^{L,T}_{r,i}\Vert &\leq \frac{1}{L^2} \sum_{j_1,j_2=1}^L \left(\Vert X(t_i,\omega_{j_1})\Vert\delta + \Vert X(t_i,\omega_{j_2})\Vert\delta + \delta^2\right)\\
    &=\delta^2 + \frac{2}{L}\sum_{j=1}^L\Vert X(t_i,\omega_j)\Vert\delta
\end{align*}
Combining these bounds leads to 
\begin{align*}
    \Vert\C^L_{r,i}-\bar\C^L_{r,i}\Vert &\leq \frac{1}{L-1}\left(L\delta^2 + \sum_{j=1}^L 2\Vert X(t_i,\omega_j)\Vert\delta \right) + \frac{L}{L-1}\left( \delta^2 + \frac{2}{L}\sum_{j=1}^L\Vert X(t_i,\omega_j)\Vert\delta \right)\\
    &= \frac{2L}{L-1}\delta^2 + \frac{4}{L-1}\sum_{j=1}^L \Vert X(t_i,\omega_j)\Vert\delta.
\end{align*}
Notice that the above expression converges as $L\to \infty$. Therefore, it is also bounded in $L$. Further, using that $t\mapsto X(t,\omega)$ is bounded almost surely and that $t\mapsto \mathbb E\Vert X(t)\Vert$ is bounded on $[0, T]$ ensures that the constants $c_1$ and $c_2$ can be chosen independent of $L$, $t_i$ and $r$.
\end{proof}
We know by Proposition \ref{pro:error_estimators} that we can control $ \Vert \rmE^L_{r,i}-\bar \rmE^L_{r,i}\Vert $ and $\Vert\C^L_{r,i}-\bar \C^L_{r,i}\Vert$ if we control the error between \eqref{eq:SDEFOM} and \eqref{eq:SDEROM}. This can, e.g., be ensured by fixing a suitable $r$. Combining this insight with Theorems \ref{thm:convergence_Expectation} and \ref{thm:convergence_Covariance}, we observe that the coefficients of \eqref{eq:SDEOIROM} computed by Algorithm \ref{alg:OpInfSDE} can be close to the matrices of the POD-ROM \eqref{eq:SDEROM} if $L$, $h$ and $r$ are chosen appropriately.

Before concluding the closeness of the ROM obtained by the presented approach to the POD-ROM, the continuity of the expectation and covariance with respect to the initial condition and the ROM coefficients needs to be established.
The version of the Gronwall lemma referred to in the proofs of the following theorems is provided in the appendix by Lemma \ref{lem:Gronwall}.
\begin{theorem}\label{thm:CovErrSol}
Let $C$ and $\hat{C}$ be the solutions of the Lyapunov differential equations
\begin{align*}
\dot{C}(t) &= \Psi(t)C(t) + C(t)\Psi(t)^T + \H(t), \quad C(0) = C_0, \\
\dot{\hat{C}}(t) &= \hat{\Psi}(t) \hat{C}(t) + \hat{C}(t)\hat{\Psi}(t)^T + \hH(t), \quad \hat{C}(0) = \hat{C}_0, 
\end{align*}
with coefficients $\H,\hH,\Psi,\hPsi\in L^1([0,T],\mathbb{R}^{n\times n})$ and let the difference of the coefficient functions be denoted by $\Delta\Psi=\hPsi-\Psi$ and $\Delta\H=\hH-\H$.
Define $\Delta C := \hat{C} - C$ and $e_{cov} = \Vert \Delta C\Vert$.
There exist constants $\alpha_c=\alpha_c(T), \beta_c=\beta_c(T),\gamma_c=\gamma_c(T)>0$ depending on the end-time $T$, $C_0$, $H$ and $\Psi$, such that 
\[e_{cov}(t) \leq \left( \alpha_c+\beta_c\int_0^t\Vert\Delta\Psi(s)\Vert\rmd s + \gamma_c\int_0^t\Vert\Delta\H(s)\Vert\rmd s \right)\exp{\left(\int_0^t\Vert\Delta\Psi(s)\Vert\rmd s\right)}\]
for all $t\in[0,T]$.
\end{theorem}
\begin{proof}
By the Carathéodory existence theorem \cite[Theorem 1.45]{roubicekNonlinearPartialDifferential2013}, the unique functions $C$ and $\hat{C}$ exist, are absolutely continuous and fulfil their respective Lyapunov differential equations for almost all $t\in[0,T]$.
From the differential equations of $C(t)$ and $\hat{C}(t)$ one obtains the evolution equation of $\Delta C(t)$
\begin{align*}
\Delta\dot{C}(t) &= \hPsi(t)\hat{C}(t)+\hat{C}(t)\hPsi(t)^T + \hH(t) - \left[\Psi(t)C(t)+C(t)\Psi(t)^T + \H(t) \right]
\end{align*}
with initial condition $ \Delta C(0) = \hat{C}_0 - C_0$ and $\Delta\dot{C}(t):= \dot{\hat{C}}(t) - \dot C(t)$. 
Substituting $\hat{C}(t) = \Delta C(t) + C(t)$, yields
\begin{align*}
\Delta\dot{C}(t)
&= \hPsi(t)(\Delta C(t) + C(t))+ (\Delta C(t) + C(t))\hPsi(t)^T + \hH(t) \\
& \quad- \left[\Psi(t)C(t)+C(t)\Psi(t)^T + \H(t) \right]\\
&= \hPsi(t)\Delta C(t)+\Delta C(t) \hPsi(t)^T + \Xi(t),
\end{align*}
with
\[\Xi(t) = \Delta\Psi(t) C(t) + C(t)\Delta\Psi(t)^T+\Delta\H(t).\]
Switching to the integral representations, $\Delta C(t)$ satisfies
\[\Delta C(t) = \Delta C(0) + \int_0^t \hPsi(s)\Delta C(s)+\Delta C(s)\hPsi(s)^T+\Xi(s)\rmd s\]
and the norm of $\Delta C(t)$ is bounded by
\[e_{cov}(t) \leq \underbrace{e_{cov}(0)+\int_0^t \Vert \Xi(s)\Vert\rmd s}_{=:\alpha(t)} + \int_0^t\underbrace{2\Vert\hPsi(s)\Vert}_{=:\beta(s)} e_{cov}(s) \rmd s.\]
The fundamental theorem of calculus for Lebesgue integrals \cite{cohnMeasureTheorySecond2013} states that $\alpha$ is absolutely continuous.
Furthermore, $e_{cov}$ is continuous, since $C$ and $\hat{C}$ are continuous.
Including the fact that the product of a continuous function and an integrable function over a compact interval is integrable again, the assumptions of the Gronwall lemma \ref{lem:Gronwall} are satisfied.
Together with the triangle inequality, one obtains explicit bounds for all $t\in[0,T]$,
\begin{align*}
 e_{cov}(t) &\leq \alpha(t)\exp{\left(\int_0^t \beta(s)\rmd s\right)}\\
 &= \left( e_{cov}(0)+\int_0^t \Vert \Xi(s)\Vert\rmd s\right)\exp{\left(\int_0^t 2\Vert\hPsi(s)\Vert\rmd s\right)}\\
 &\leq \left( e_{cov}(0) + \int_0^t 2\Vert \Delta\Psi(s)\Vert\Vert C(s)\Vert + \Vert\Delta\H(s)\Vert\rmd s\right)\\
 &\quad\quad\cdot\exp{\left(\int_0^t \Vert\Psi(s)\Vert + \Vert\Delta\Psi(s)\Vert\right) \rmd s },
\end{align*}
which becomes 
\begin{equation*}
 e_{cov}(t) \leq \left(\alpha_c+\tilde\beta_c\int_0^t\Vert\Delta\Psi(s)\Vert\rmd s + \gamma_c\int_0^t\Vert\Delta\H(s)\Vert\rmd s \right)\exp{\left(\int_0^t\Vert\Delta\Psi(s)\Vert\rmd s\right)}
\end{equation*}
for all $t\in[0,T]$, with 
\begin{align*}
 \alpha_c = \Vert \hat C(0)-C(0)\Vert\gamma_c, \quad 
 \tilde\beta_c = 2 \tilde c_1 \gamma_c, \quad
 \gamma_c = e^{\int_0^T\Vert \Psi(s)\Vert \rmd s}
\end{align*}
and $\tilde c_1 = \max_{s\in[0,T]}\Vert C(s)\Vert$.
By bounding $\Vert C(t)\Vert$, one can remove the dependence of $\tilde\beta_c$ on the  state $C(t)$ and obtain constants that depend solely on the coefficient functions and the difference in initial values.
Employing the Gronwall lemma \ref{lem:Gronwall} on the norm of 
\[C(t) = C_0 + \int_0^t \Psi(s)C(s)+C(s)\Psi(s)^T + \H(s)\rmd s\]
yields
\begin{align*}
 \Vert C(t) \Vert \leq \left(\Vert C_0\Vert +\int_0^T \Vert \H(s)\Vert \rmd s \right)\exp{\left( 2\int_0^T \Vert\Psi(s)\Vert\rmd s\right)} =: c_1 
\end{align*}
for all $0\leq t\leq T$.
Defining the constant $\beta_c = 2c_1 \gamma_c$ then completes the proof.
\end{proof} 
\begin{theorem}\label{thm:ExpErrSol}
Let $E$ and $\hat{E}$ be the solutions of the ODEs
\begin{align*}
 \dot{E}(t) &= \A(t)E(t)+\B(t),\quad E(0) = E_0,\\
 \dot{\hat{E}}(t) &= \hA(t)\hat{E}(t)+\hB(t),\quad \hat{E}(0) = \hat{E}_0,
\end{align*}
with coefficients $\A,\hA\in L^1([0,T],\Rnn)$ and $\B,\hB\in L^1([0,T],\Rn)$.
Let the deviations in the coefficients be defined by $\Delta \A = \hA-\A$ and $\Delta\B = \hB-\B$ and the difference between the solutions by $\Delta E =\hat{E}-E$ and $e_{exp} = \Vert \Delta E\Vert$.
There exist constants $\alpha_e,\beta_e,\gamma_e\geq 0$ depending on the end-time $T$, $A$ and $B$, such that \[ e_{exp}(t) \leq \left( \alpha_e+\beta_e\int_0^t\Vert\Delta\A(s)\Vert\rmd s + \gamma_e\int_0^t\Vert\Delta\B(s)\Vert\rmd s \right)\exp{\left(\int_0^t\Vert\Delta\A(s)\Vert\rmd s\right)}\]
for all $t\in[0,T]$.
\end{theorem}
\begin{proof}
The proof follows similar arguments as in Theorem \ref{thm:CovErrSol}.
The constants are given by
\begin{align*}
        \alpha_e = \Vert \hat{E}(0)-E(0)\Vert\gamma_e,\quad 
        \beta_e = c_3\gamma_e, \text{ and }
        \gamma_e = e^{\int_0^T\Vert\A(s)\Vert\rmd s},
    \end{align*}
    where 
\[c_3 =\left( E_0 + \int_0^T \Vert \B(s) \Vert \rmd s\right)\exp{\left( \int_0^T \Vert \A(s)\Vert \rmd s\right)}.\]
\end{proof}

By Theorems \ref{thm:CovErrSol} and \ref{thm:ExpErrSol} the distance of two systems of the structure of \eqref{eq:SDEROM} in the first and second moment are determined by the distance of their respective system coefficients.
In our case, the coefficients are
\begin{align*}
    \Psi (t) = \A_r + \sum_{i=1}^m \N_{r,i}u_i(t), \quad 
    \hPsi (t) = \hA_r^{L,h} + \sum_{i=1}^m \hN_{r,i}^{L,h}u_i(t),\quad
    \H(t) = \H_r,\quad
    \hH(t) = \hH_r^{L,h}
\end{align*}
and 
\begin{align*}
    \A (t) = \Psi(t),\quad
    \hA (t) = \hPsi(t),\quad 
    \B(t) = \B_r u(t)\quad
    \hB(t) = \hB^{L,h}_r u(t).
\end{align*}
By Theorems \ref{thm:convergence_Expectation}, \ref{thm:convergence_Covariance} and Proposition \ref{pro:error_estimators} we know that, given an error threshold $\epsilon$, there exists a $L=L_{\epsilon}, h=h_{\epsilon}$ and a ROM dimension $r=r_{\epsilon}$, such that the distance of inferred coefficients to the intrusive coefficients is less than $\epsilon$ with probability $1$.
Thus, for instance,
\begin{align*}
    \int_0^T\Vert \Psi(s)-\hPsi(s)\Vert\rmd s &= \int_0^T\Vert \A_r-\hA_r^{L,h} + \sum_{i=1}^m \left( \N_{r,i}-\hN_{r,i}^{L,h}\right)u_i(s)\Vert\rmd s\\
    &\leq \Vert \A_r-\hA_r^{L,h}\Vert T + \sum_{i=1}^m \Vert \N_{r,i}-\hN_{r,i}^{L,h}\Vert \int_0^T \Vert u_i(s)\Vert\rmd s\\
    &\leq c\epsilon,
\end{align*}
for some constant $c=c(T,u)$ with probability $1$.
For sufficiently large ROM dimension $r$ and sample size $L$, as well as a sufficiently small time step size $h$, it can therefore be guaranteed that the expectation and covariance of the OpInf ROM can be arbitrarily close to the expectation and covariance of the POD-ROM.
Including the fact that the ROMs are Gaussian for each fixed time $t\in[0,T]$, this means that, under the appropriate conditions, the presented OpInf approach for linear SDEs with additive Gaussian noise produces a ROM that is close in distribution to the intrusive POD-ROM with probability $1$.

\section{Numerical experiments}
\label{sec:Experiments}
Due to the slow rate of convergence of the empirical estimators, we can expect a large number of FOM samples to be required to approximate the expectation and covariance well.
While this is computationally expensive, especially if the FOM dimension is very large, constructing a ROM using the presented method is still beneficial. 
For instance, in optimal control problems, where we need to evaluate the system for many controls. Now, if we approximate the expectation of a quantity related to the state for each of the many controls, we have a huge sampling effort, which is the number of samples in a Monte-Carlo approach times the number of controls we are interested in. This is usually not feasible in the FOM setting. 
Rather than using the computational resources for repeated FOM evaluations within an iterative solver, we can leverage these samples to construct the OpInf ROM using a small set of inputs combined with a moderately large number of samples (smaller than in the above mentioned Monte-Carlo method). The inferred ROM can then serve as a surrogate model in the initial iterations of the optimisation algorithm, reducing computational costs while being close in distribution to the FOM.

In the following section, two experiments are presented, with the help of which the developed method is compared against POD.
The drift and diffusion coefficients of the FOM are obtained by appropriate finite difference discretisations of the specified PDEs in the spatial coordinates.
Though strictly speaking, the time derivative of a Wiener process does not exist, the diffusion coefficient $\M\rmd W(t)$ of the FOM can be thought of as corresponding to a term $\sigma(x) \dot{W}(t)$.
This notation is used for the remainder of this section for ease of notation and brevity.
We choose the number of time steps $s$ and their length $h$ such that the end time point $T=sh=1$ in both examples.

\paragraph{1d Heat equation}
The drift coefficients of the test-model are obtained by spatial discretisation of the one dimensional heat equation with Dirichlet boundary conditions
\begin{equation} \label{eq:FOMHeat}
\begin{split}
    \frac{\partial}{\partial t} y(x, t) &= 0.1\frac{\partial^2}{\partial^2 x}y(x, t)+u(t)\frac{\partial}{\partial x}y(x, t)+\sigma(x)\dot W(t),\quad x\in(0,1),t\in(0,T],\\
y(x,0) &= y_0(x) \quad \text{ and }\quad y(0,t) = y(1,t) = u(t),
\end{split}
\end{equation}
by finite differences using $n=100$ points.
In this example, the noise generating Wiener process is of dimension $d=2$.
The columns of $\M$ correspond to the evaluation of $\sigma_1(x)=0.1\exp (-10(x-\tfrac{1}{2})^2)$ and $\sigma_2(x) = 0.1\sin(2\pi x)$ on the grid $x_1,\dots,x_n$.
In this experiment, $s=1000$ equidistant time steps of length $h=10^{-3}$ were used.

\paragraph{2d Heat equation}
This example models the heat spread on the unit square with a hole.
The input, modelling the source temperature, is assumed to be noisy and, as such, is modelled by a noisy control $\tilde{u}(t) = u(t)+\dot{W}(t)$.
Here, $u$ is a deterministic function of time which is perturbed by a $d=1$ dimensional Gaussian noise at each time-point $t$.
The spatial domain, on which the dynamics are defined, is non-convex and the Dirichlet boundary condition is non-continuous at two points of the boundary. 
The mathematical formulation is given by
\begin{equation}\label{eq:FOM2dHeat}
    \begin{aligned}
    \mathbbm{1}_{\Omega_1}(x)\tilde u(t)&=\frac{\partial }{\partial t}y(x,t) -  0.01\Delta y(x,t), \quad &\forall (x,t) \in (\Omega\setminus\Omega_2)\times [0,T]\\
    \mathbbm{1}_{\partial\Omega_1}(x)\tilde u(t)&=y(x,t) &\forall (x,t)\in \partial\Omega\times [0,T]\\ 
    0&=y(x,t) &\forall (x,t)\in\Omega_2\times[0,T]\\
    0&=y(x,0) &\forall x \in \Omega,
    \end{aligned}
\end{equation} 
\newpage
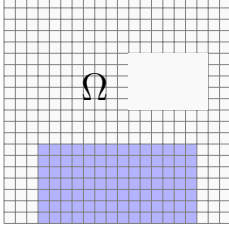
\begin{wrapfigure}{l}{0.3\textwidth}
\begin{tikzpicture}[scale=3] 
  \fill[gray!5] (0,0) rectangle (1,1);
  \fill[blue!30] (0.15,0) rectangle (0.85,0.35);
  \draw[step=0.05,gray,very thin] (0,0) grid (1,1); 
  \fill[gray!5] (0.55,0.5) rectangle (0.9,0.75);
  \node at (0.4,0.6) {\Large$\Omega$};
\end{tikzpicture}
\end{wrapfigure}
\noindent where the respective domains are specified by 
$\Omega = [0,1]^2$, $\Omega_1=[0.12,0.88]\times[0,0.36], \Omega_2 = [0.51,0.48]\times[0.94,0.79]$ and $\mathbbm{1}_A(x)$ is the indicator function of a set $A\subset\Omega$ evaluated at $x\in\Omega$.
The domain $\Omega\setminus\Omega_2$ was discretised using $n=1016$ equidistant points.
Due to the construction of the example, the matrices $\B\in\Rn$ and $\M\in\Rn$ are identical and consist of only one column.
In this setup, the PDE dynamics evolve in a non-convex domain. 
The Dirichlet boundary condition is piece-wise continuous on $\partial\Omega$ with non-continuities at $x\in\{(0.12,0),(0.88,0)\}$.
In this experiment, $s=100$ equidistant time steps of length $h=10^{-2}$ were used.

\begin{figure}
    \centering
    \setlength\fwidth{0.4\textwidth}
    \setlength\fheight{2cm}
    \begin{subfigure}[t]{0.5\textwidth}
       \begin{tikzpicture}
\begin{axis}[%
width=0.951\fwidth,
height=\fheight,
at={(0\fwidth,0\fheight)},
scale only axis,
xmin=0,
xmax=1,
separate axis lines,
every outer y axis line/.append style={black},
every y tick label/.append style={font=\color{black}},
every y tick/.append style={black},
ymin=-0.2,
ymax=0.2,
yminorticks=true,
axis background/.style={fill=white},
yticklabel pos=left,
xmajorgrids,
ymajorgrids,
ylabel= {$i\equal 1$},
title = left singular vectors of $\bbX$ and $\F^L$,
legend style={legend cell align=left, align=left, draw=white!15!black},
legend pos = north west,
]
\addplot [color=red, line width=2.0pt] table [y=U1, x=x] {LSV_X.dat}; 
\addplot [color=black, dashed, line width=2.0pt] table [y=U1, x=x] {LSV_F.dat}; 
\legend{$v_i$,$v_{m,i}$}
\end{axis}
\end{tikzpicture}
    \end{subfigure}%
    \begin{subfigure}[t]{0.5\textwidth}
        \begin{tikzpicture}
\begin{axis}[%
width=0.951\fwidth,
height=\fheight,
at={(0\fwidth,0\fheight)},
scale only axis,
xmin=0,
xmax=1,
separate axis lines,
every outer y axis line/.append style={black},
every y tick label/.append style={font=\color{black}},
every y tick/.append style={black},
ymode=log,
ymin=1e-07,
ymax=1e-2,
yminorticks=true,
axis background/.style={fill=white},
yticklabel pos=left,
xmajorgrids,
ymajorgrids,
title = pointwise difference
]
\addplot [color=blue, dotted, line width=2.0pt] table [y=U1, x=x] {LSV_diff.dat}; 
\end{axis}
\end{tikzpicture}
    \end{subfigure}%
   \vspace{-0.5cm}
    \hfill
    \begin{subfigure}[t]{0.5\textwidth}
       \begin{tikzpicture}
\begin{axis}[%
width=0.951\fwidth,
height=\fheight,
at={(0\fwidth,0\fheight)},
scale only axis,
xmin=0,
xmax=1,
separate axis lines,
every outer y axis line/.append style={black},
every y tick label/.append style={font=\color{black}},
every y tick/.append style={black},
ymin=-0.2,
ymax=0.2,
yminorticks=true,
axis background/.style={fill=white},
yticklabel pos=left,
xmajorgrids,
ymajorgrids,
ylabel= {$i\equal 2$},
legend style={legend cell align=left, align=left, draw=white!15!black}
]
\addplot [color=red, line width=2.0pt] table [y=U2, x=x] {LSV_X.dat}; 
\addplot [color=black, dashed, line width=2.0pt] table [y=U2, x=x] {LSV_F.dat}; 
\end{axis}
\end{tikzpicture}
    \end{subfigure}%
    \begin{subfigure}[t]{0.5\textwidth}
        \begin{tikzpicture}
\begin{axis}[%
width=0.951\fwidth,
height=\fheight,
at={(0\fwidth,0\fheight)},
scale only axis,
xmin=0,
xmax=1,
separate axis lines,
every outer y axis line/.append style={black},
every y tick label/.append style={font=\color{black}},
every y tick/.append style={black},
ymode=log,
ymin=1e-07,
ymax=1e-2,
yminorticks=true,
axis background/.style={fill=white},
yticklabel pos=left,
xmajorgrids,
ymajorgrids,
]
\addplot [color=blue, dotted, line width=2.0pt] table [y=U2, x=x] {LSV_diff.dat}; 
\end{axis}
\end{tikzpicture}
    \end{subfigure}%
   \vspace{-0.5cm}
    \hfill
     \begin{subfigure}[t]{0.5\textwidth}
       \begin{tikzpicture}
\begin{axis}[%
width=0.951\fwidth,
height=\fheight,
at={(0\fwidth,0\fheight)},
scale only axis,
xmin=0,
xmax=1,
separate axis lines,
every outer y axis line/.append style={black},
every y tick label/.append style={font=\color{black}},
every y tick/.append style={black},
ymin=-0.2,
ymax=0.2,
yminorticks=true,
axis background/.style={fill=white},
yticklabel pos=left,
xmajorgrids,
ymajorgrids,
ylabel= {$i\equal 10$},
legend style={legend cell align=left, align=left, draw=white!15!black}
]
\addplot [color=red, line width=2.0pt] table [y=U10, x=x] {LSV_X.dat}; 
\addplot [color=black, dashed, line width=2.0pt] table [y=U10, x=x] {LSV_F.dat}; 
\end{axis}
\end{tikzpicture}
    \end{subfigure}%
    \begin{subfigure}[t]{0.5\textwidth}
        \begin{tikzpicture}
\begin{axis}[%
width=0.951\fwidth,
height=\fheight,
at={(0\fwidth,0\fheight)},
scale only axis,
xmin=0,
xmax=1,
separate axis lines,
every outer y axis line/.append style={black},
every y tick label/.append style={font=\color{black}},
every y tick/.append style={black},
ymode=log,
ymin=1e-07,
ymax=1e-2,
yminorticks=true,
axis background/.style={fill=white},
yticklabel pos=left,
xmajorgrids,
ymajorgrids,
]
\addplot [color=blue, dotted, line width=2.0pt] table [y=U3, x=x] {LSV_diff.dat}; 
\end{axis}
\end{tikzpicture}
    \end{subfigure}%
    \caption{
    Comparison of left singular vectors of $\bbX$ and $\F^L$ for the 1d Heat example with $n =100$, a zero initial condition and $u(t)=\cos (\frac{2\pi t}{T}),t\in[0,T]$.
    The data was generated using $L=10^4, s=1000$ and $h=10^{-3}$.
    The vectors $v_i$ and $v_{m,i}$, $i=1,2,10$ are shown in the left column. 
    The right column shows the pointwise absolute difference between the respective vectors.}
    \label{fig:SingularVectors}
\end{figure}
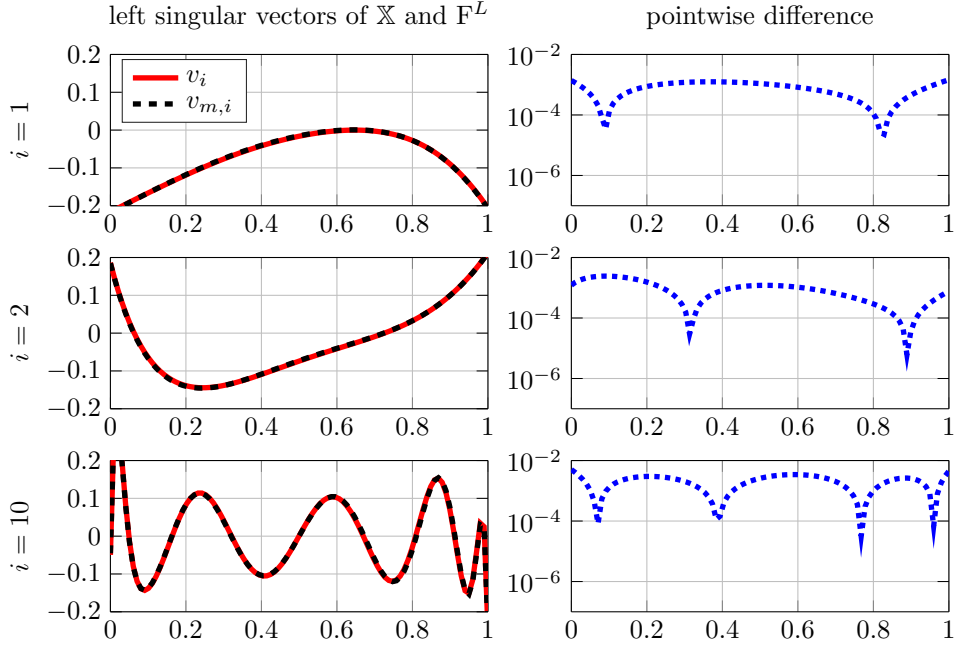

The FOMs obtained from the described experiments are of the form
\begin{equation}\label{eq:FOMHeatSDE}
    \rmd X(t) = \left[ \A X(t) + \B u(t) + \N X(t)u(t)\right]\rmd t + \M \rmd W(t),\quad X(0) = X_0,
\end{equation}
where $X(t)\in\R^n, \A,\N\in\R^{n\times n}, \B\in\R^{n\times m}, \M\in \R^{n\times d},$ and $W(t)\in\R^d$ is a $d$ dimensional Wiener process.
In the case of the $2d$ Heat equation example, the coefficient $\N$ is zero.
In all examples, the input functions are $m=1$ dimensional and the correlation matrices of the noise-generating Wiener processes are the identity $\I_d$.
To perform the numerical integration, a strong drift-implicit Euler-Maruyama scheme \cite[Chapter 12.2]{kloeden_numerical_1992} with $s$ steps of size $h$ is used, which means that the final time-point in both cases is $T=sh=1$.
Except for the specified parameters, the respective setups are identical for both experiments.
\paragraph{Subspace computation}
For the comparison with POD in Figures \ref{fig:1dHeat} and \ref{fig:2dHeat}, the subspace $\V_r$ was constructed from the left singular vectors of the state snapshot matrix $\bbX$.
To this end, $L=10^4$ snapshots were generated using a zero initial condition, the input $u(t) = \cos(\tfrac{2\pi}{T}t), t\in[0,T],$ and, for each experiment, the respective number of time steps $s$ of length $h$.
Figure \ref{fig:SingularVectors} displays the left singular vectors $v_i$ and $v_{m,i}$ of $\bbX$ and $\F^L$ in the case of the bilinear 1d Heat equation, respectively. 
The decay rates of the respective singular values do not differ substantially.
While the plots of Figure \ref{fig:SingularVectors} in this specific case suggest that the left singular vectors obtained from $\bbX $ and $\F^L$ are (almost) identical, this is generally not the case. 
In additional experiments that are not part of this paper, we observed that the left singular vectors of $\F^L$ coincide with those of $\bbX$, up to a permutation.

\paragraph{Data generation}
In order to construct a data-matrix of small condition number, we repeatedly evaluate the FOM using $s=1000,h=10^{-3}$ for the 1d Heat example and $s=100,h=10^{-2}$ for the 2d Heat example over pairs of initial conditions and inputs.
In the first stage the FOM is queried using a zero initial condition and $k=21$ input functions $u^{(i)},\,i=1,\dots,k$.
The inputs are chosen as constant functions $u^{(i)} \equiv -2 + \tfrac{4}{k}i$.
The second part consists of querying the FOM $r$ times, each with one of the $r$ left singular vectors $v_1,\dots, v_r$ as the initial condition and a zero control $u^{(0)}$. 
For each of these $r+k$ pairs, $L=10^1, 10^2$, and $10^4$ samples are drawn. 
As the FOM state dimension is typically very large, storing the sampled trajectories or their estimated moments might be infeasible. 
To address this, we use the previously computed matrix $\V_r=[v_1,\dots,v_r]\in\R^{n\times r}$ to first project the FOM observations and then estimate and store the expectation and covariance of these projected trajectories.
The sampling algorithm is implemented in such a way that, for a given pair, the first $10$ realisations of the noise generating Wiener process are identical across all runs and the first $100$ realisations are shared between the $L=10^2$ and $L=10^4$ runs.
The realisations are independent between different pairs of initial condition and control.
Using a different number of samples to estimate the empiric moments did not effect the condition number of the data-matrices significantly.
The approximation of the time derivatives is carried out using a central finite difference scheme.

\paragraph{Testing setup}
The quality of the obtained ROMs is compared using the relative error in expectation, covariance and in the weak sense. 
Concretely, this means that the following quantities 
\begin{subequations}
\label{eq:DefinitionErrors}
\begin{align}
    e_E &=\frac{\sum_{i=1}^s\| E(t_i)-\V_rE_r(t_i)\|_2^2}{\sum_{i=1}^s\| E(t_i)\|_2^2},\\  
e_C &=\frac{\sum_{i=1}^s\| C(t_i)-\V_rC_r(t_i)\V_r^T\|_F^2}{\sum_{i=1}^s\| C(t_i)\|_F^2}, \\
e_{\phi,i}(\tau) &= \frac{\|\E{\phi_i(X(\tau))}-\E{\phi_i(\V_rX_r(\tau))}\| }{\E{\phi_i (X(\tau))}} \quad \text{ for some } 0\leq \tau\leq T,
\end{align}
\end{subequations}
were computed for the POD-ROM and the ROM obtained by OpInf. 
The third quantity is called the \emph{relative weak error} of $X$ with respect to functionals $\Phi_i:\Rn\to \R$, cf. \cite{kloeden_numerical_1992}.
In the presented error plots of Figures \ref{fig:1dHeat} and \ref{fig:2dHeat}, the weak error is measured at the end-time $\tau=T$.
The functionals 
\begin{equation}\label{eq:ErrorFunctionals}
\phi_1 (\tau) = \Vert X(\tau)\Vert_2^2, \text{ and }
\phi_2 (\tau) = \frac{1}{n}\sum_{i=1}^n X_i(\tau)^3e^{X_i(\tau)}
\end{equation}
were chosen for the computation of the weak error.
While $\phi_1$ refers to the second moment of $X(t)$, the influence of $X$ on values of $\phi_2$ is dominated by the large components of $X$.
The number of samples used for the error computation is set to $10^6$.
The ROMs are tested against the FOM for a zero initial condition and the control function $u(t) = \cos(\tfrac{5\pi}{T}t),\,t\in[0,T]$.
The random number generator seed was fixed and reset before each estimation of the empiric moments of the FOM and the ROMs.
Consequently, the FOM and POD-ROM are driven by the same noise realisations. 
This is also the case for the OpInf ROM, as long as the inferred noise dimension $d_r$ is equal to $d$. 
In our experiments, the only case where $d_r\neq d$ was for the OpInf ROM of dimension $r=1$ in the one dimensional Heat example.
The seed used for testing was set to an arbitrary value different from the seed used generating the data.

\begin{figure}
    \setlength\fwidth{0.4\textwidth}
    \setlength\fheight{4cm}
    \begin{subfigure}[b]{0.5\textwidth}
        \begin{tikzpicture}
            \begin{axis}[ymin=1e-2,ymax=1000,ymode=log,width=\fwidth,height=\fheight,scale only axis, xlabel=time step $i$,title= $1d$ Heat equation \eqref{eq:FOMHeat}, ylabel style={yshift=-1em}]
                \addplot[dotted, line width = 1pt] table[x=x, y=normE] {data_edit2/snr_1dHeat.txt};
                \addplot[dashed,line width = 1pt] table[x=x, y=normC] {data_edit2/snr_1dHeat.txt};
                \addplot[solid,line width = 1pt] table[x=x, y=snr] {data_edit2/snr_1dHeat.txt};
            \end{axis}
        \end{tikzpicture} 
    \end{subfigure}%
    \hspace{0.2cm}
    \begin{subfigure}[b]{0.5\textwidth}
        \begin{tikzpicture}
            \begin{axis}[ymin=1e-2,ymax=1000,ymode=log,width=\fwidth, height=\fheight,scale only axis,xlabel=time step $i$,title= $2d$ Heat equation \eqref{eq:FOM2dHeat}, 
            legend style={fill=none,draw=none}]
                \addplot[dotted, line width = 1pt] table[x=x, y=normE] {data_edit2/snr_2dHeat.txt};
                \addplot[dashed, line width = 1pt] table[x=x, y=normC] {data_edit2/snr_2dHeat.txt};
                \addplot[solid, line width = 1pt] table[x=x, y=snr] {data_edit2/snr_2dHeat.txt};
                \legend{$\Vert \rmE^L_i\Vert_2$,$\Vert\C^L_i\Vert_F$,$\mathrm{SNR}(t_i)$}
            \end{axis}
        \end{tikzpicture} 
    \end{subfigure}%
    \hfill
    \caption{ Norms of expectation and covariance as well as the signal-to-noise ration $\mathrm{SNR}(t)$ of FOM state snapshots used for subspace construction for the $1d$ and $2d$ Heat equation, respectively.}
    \label{fig:signal_to_noise}
\end{figure}
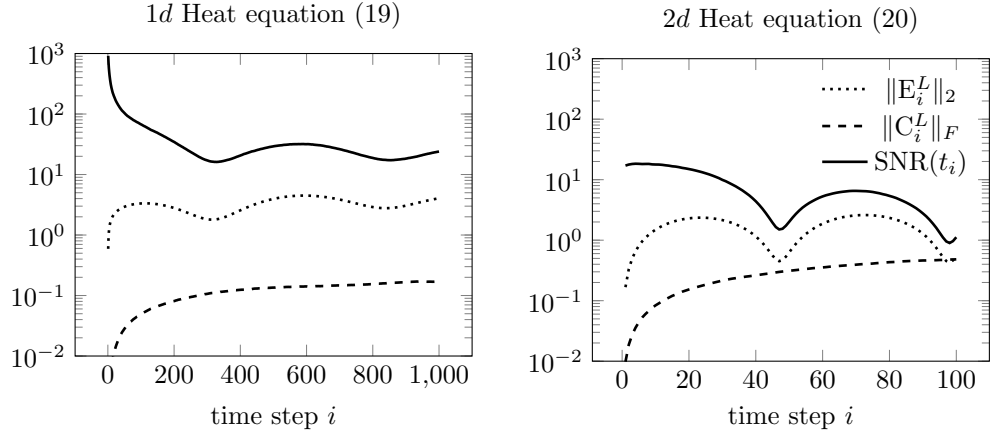

\paragraph{Results}
Figures \ref{fig:1dHeat} and \ref{fig:2dHeat} display the relative errors in the moments and relative weak errors for the $1d$ and $2d$ Heat equation, respectively.
The estimated signal-to-noise ratio \[\mathrm{SNR}(t_i)=\frac{\Vert\rmE_{f,i}^L\Vert_2}{\Vert\C_{f,i}^L\Vert_F}\] as well as the 2-norm of the expectation $\rmE^L_i$ and the Frobenius norm of $\C^L_i$ of the data used for the subspace computation are displayed in Figure \ref{fig:signal_to_noise}.
As the initial condition and control are the same for generating the state snapshots and computing the test set, this ratio describes the relative importance of expectation and covariance in the test set as well.
Figure \ref{fig:signal_to_noise} reveals that, compared to the $1d$ Heat equation \eqref{eq:FOMHeat}, the contribution of the noise in the data relative to the expectation is between one and two orders of magnitude larger in the $2d$ Heat example \eqref{eq:FOM2dHeat}.
In contrast, the norms of the covariances are similar between both examples, with $\smash{\Vert \C^L_i}\Vert_F$ exhibiting a slightly larger increase.
An immediate observation of Figure \ref{fig:1dHeat} is that for the $1d$ Heat equation the ROM approximation error of the covariance only decreases significantly after $r=4$.
The same can be observed for the ROM dimensions $r=1,2$ in the $2d$ Heat example.
As this behaviour is present in all of the OpInf ROMs as well as the POD-ROM, this suggests that the that the column space of $\bbX$ that corresponds to the FOM covariance is less significant compared to the column space associated to the FOM expectation in both examples, with the difference being smaller for the $2d$ Heat eq.
This is supported by computing the expectation $\rmE^L$ and covariance $\C^L$ of the data contained in the state snapshot matrix $\bbX$.
We find that, while for the first example $\Vert \rmE^L\Vert_F \approx 103$ and $\Vert\C^L\Vert_F \approx 4$, for the $2d$ Heat example $\Vert\rmE^L\Vert_F\approx 17 $ and $\Vert\C^L\Vert_F\approx 3$, which agrees with the previous observations of Figure \ref{fig:signal_to_noise}.
In both experiments the OpInf method presented in this paper provides better approximations in the expectation and covariance as the number of samples increases.
Generally, it can be observed that the OpInf ROMs obtained from $L=100$ and $L=10000$ samples provide similar errors for the $1d$ Heat equation and for the expectation approximation error of the $2d$ Heat equation.
This closeness indicates to us that for these quantities either only $100$ samples are necessary to reduce the sampling error in the data below the thresholds arising in the approximated deterministic inference problems, or that drawing $10^6$ samples to estimate the errors does not provide a fine enough resolution.
Executing Algorithm \ref{alg:OpInfSDE} again on the respective FOMs with the diffusion coefficients set to zero reproduces the observed errors for $L=10000$ in both cases --- except for the covariance error metric, which cannot be tested using this method.
This means, that regarding the error metrics except the covariance approximation error, this ROM achieves the same error as if we would have been able to observe the FOM expectation directly.
Consequently, for these two examples, it is sufficient to use 100 samples to achieve a good approximation of the expected value.
In contrast, the covariance approximation error of the second example illustrates that even though there is no reduction in the expectation error when increasing $L$ from 100 to 10000, the covariance approximation error might still decrease.
Interestingly, the OpInf ROM obtained from $L=10$ samples provides relatively good approximation results, despite using one and three orders of magnitude fewer samples, respectively.
However, this particular observation is highly dependent on the seed chosen during the data generation.
Conversely, as it is to be expected, the robustness of the inference to the noise realisation increases for larger sample sizes. 
\begin{figure}[h!]
 \centering
    \newcommand{\xmin}{1}
    \newcommand{\xmax}{10}
    \newcommand{\ymin}{1e-4}
    \newcommand{\ymax}{1e0}
    \setlength\fwidth{0.4\textwidth}
    \setlength\fheight{5cm}
    \begin{subfigure}[b]{0.5\textwidth}
        \begin{tikzpicture}
\begin{axis}[%
width=0.951\fwidth,
height=\fheight,
at={(0\fwidth,0\fheight)},
scale only axis,
unbounded coords=jump,
title = error in expectation $e_E$,
legend pos=south west,
legend style={fill=none, draw=none, text width=3cm, align=left,yshift = -0.2cm, xshift= -0.1cm},
xmin=\xmin,
xmax=\xmax,
ymode=log,
ymin=\ymin,
ymax=\ymax,
xlabel = ROM dimension $r$,
yminorticks=true,
axis background/.style={fill=white},
ymajorgrids=true, 
xmajorgrids=true
]
    \addplot[mark=square*, cPOD] table [y=err1, x=rank] {data_edit2/Heat1d/E_error_POD.txt};     
    \addplot[mark=*, cOI1] table [y=err1, x=rank] {data_edit2/Heat1d/E_error_OpInf.txt}; 
    \addplot[mark=*, cOI2] table [y=err2, x=rank] {data_edit2/Heat1d/E_error_OpInf.txt};
    \addplot[mark=*, cOI3] table [y=err3, x=rank] {data_edit2/Heat1d/E_error_OpInf.txt};
\end{axis}
\end{tikzpicture}%
    \end{subfigure}%
    \begin{subfigure}[b]{0.5\textwidth}
        \begin{tikzpicture}
\begin{axis}[%
width=0.951\fwidth,
height=\fheight,
at={(0\fwidth,0\fheight)},
scale only axis,
unbounded coords=jump,
title = error in covariance $e_C$,
legend pos=south west,
legend style={
fill=none,
draw=none, 
text width=2.3cm, 
font = \small,
align=left,
},
xmin=\xmin,
xmax=\xmax,
ymode=log,
ymin=\ymin,
ymax=\ymax,
xlabel = ROM dimension $r$,
yminorticks=true,
axis background/.style={fill=white},
ymajorgrids=true, 
xmajorgrids=true
]
    \addplot[mark=square*, cPOD] table [y=err1, x=rank] {data_edit2/Heat1d/C_error_POD.txt};     
    \addplot[mark=*, cOI1] table [y=err1, x=rank] {data_edit2/Heat1d/C_error_OpInf.txt}; 
    \addplot[mark=*, cOI2] table [y=err2, x=rank] {data_edit2/Heat1d/C_error_OpInf.txt};
    \addplot[mark=*, cOI3] table [y=err3, x=rank] {data_edit2/Heat1d/C_error_OpInf.txt};
    \legend{POD,OpInf $L=10^1$,OpInf $L=10^2$, OpInf $L=10^4$}
\end{axis}
\end{tikzpicture}%
    \end{subfigure}%
    \hfill
    \vspace{0.2cm}
    \renewcommand{\ymin}{1e-6}
    \begin{subfigure}[b]{0.5\textwidth}
        \begin{tikzpicture}
\begin{axis}[%
width=0.951\fwidth,
height=\fheight,
at={(0\fwidth,0\fheight)},
scale only axis,
unbounded coords=jump,
title = weak error $e_{\phi,1}(T)$,
xmin=\xmin,
xmax=\xmax,
ymode=log,
ymin=\ymin,
xlabel = ROM dimension $r$,
ymax=\ymax,
yminorticks=true,
axis background/.style={fill=white},
ymajorgrids=true, 
xmajorgrids=true,
]
    \addplot[mark=square*, cPOD] table [y=err1, x=rank] {data_edit2/Heat1d/f1_error_POD.txt};     
    \addplot[mark=*, cOI1] table [y=err1, x=rank] {data_edit2/Heat1d/f1_error_OpInf.txt}; 
    \addplot[mark=*, cOI2] table [y=err2, x=rank] {data_edit2/Heat1d/f1_error_OpInf.txt};
    \addplot[mark=*, cOI3] table [y=err3, x=rank] {data_edit2/Heat1d/f1_error_OpInf.txt};
\end{axis}
\end{tikzpicture}%
    \end{subfigure}%
    \begin{subfigure}[b]{0.5\textwidth}
        \begin{tikzpicture}
\begin{axis}[%
width=0.951\fwidth,
height=\fheight,
at={(0\fwidth,0\fheight)},
scale only axis,
unbounded coords=jump,
legend pos=south west,
legend style={fill=none, draw=none, text width=3cm, align=left,yshift = -0.2cm, xshift= -0.1cm},
title = weak error $e_{\phi,2}(T)$,
xmin=\xmin,
xmax=\xmax,
ymode=log,
ymin=\ymin,
xlabel = ROM dimension $r$,
ymax=\ymax,
yminorticks=true,
axis background/.style={fill=white},
ymajorgrids=true, 
xmajorgrids=true
]
    \addplot[mark=square*, cPOD] table [y=err1, x=rank] {data_edit2/Heat1d/f2_error_POD.txt};     
    \addplot[mark=*, cOI1] table [y=err1, x=rank] {data_edit2/Heat1d/f2_error_OpInf.txt}; 
    \addplot[mark=*, cOI2] table [y=err2, x=rank] {data_edit2/Heat1d/f2_error_OpInf.txt};
    \addplot[mark=*, cOI3] table [y=err3, x=rank] {data_edit2/Heat1d/f2_error_OpInf.txt};
\end{axis}
\end{tikzpicture}%
    \end{subfigure}%
    \caption{$1d$ Heat equation. The subspace computation, data generation and testing were performed with $s=1000$ time steps of size $h=10^{-3}$. The number of samples drawn to estimate the error quantities was $10^6$.}
    \label{fig:1dHeat}
\end{figure}
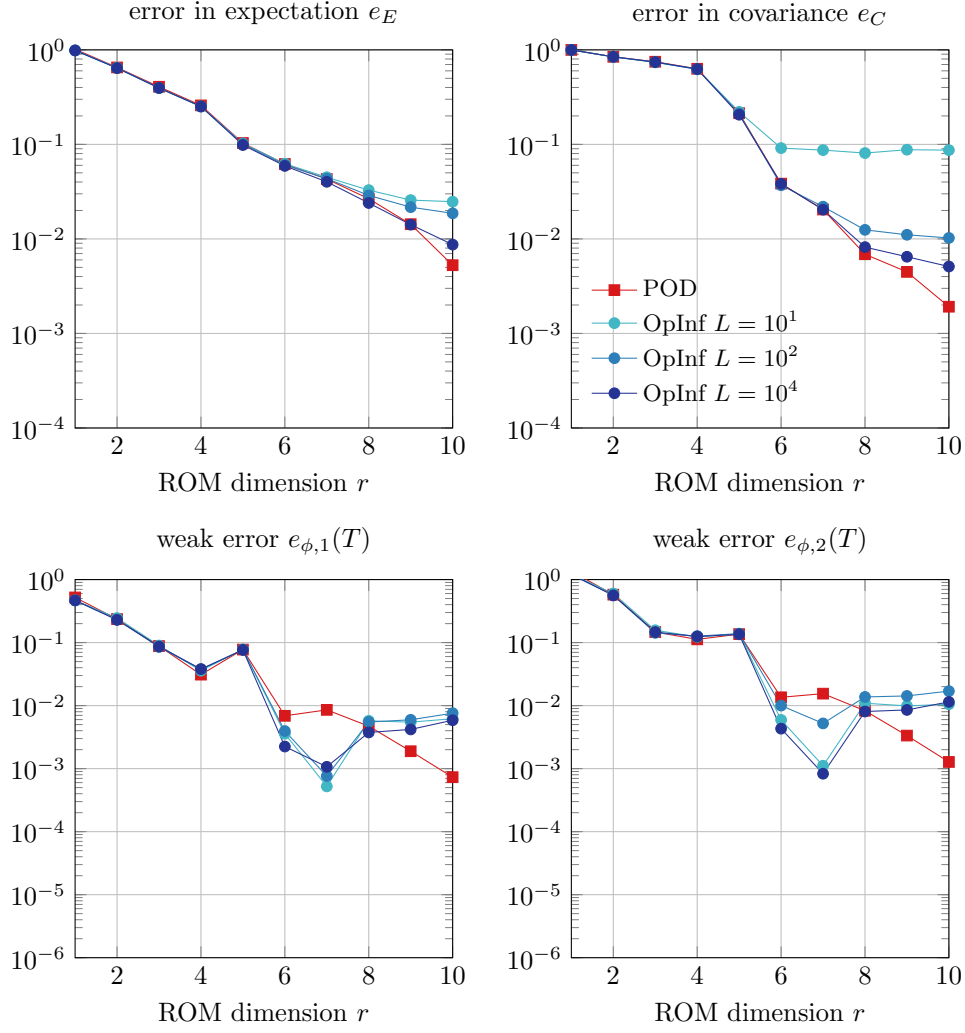
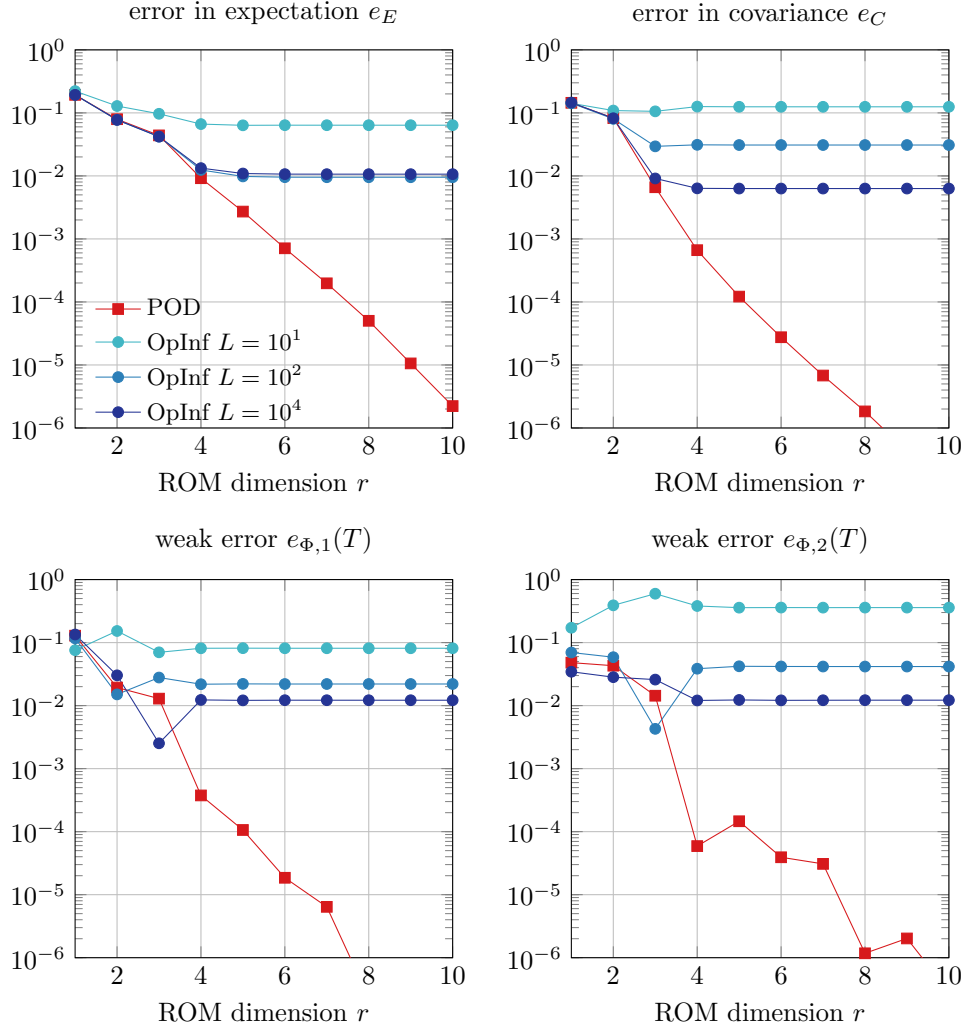
\begin{figure}[h!]
    \centering
    \newcommand{\xmin}{1}
    \newcommand{\xmax}{10}
    \newcommand{\ymin}{1e-8}
    \newcommand{\ymax}{1e0}
    \setlength\fwidth{0.4\textwidth}
    \setlength\fheight{5cm}
    \renewcommand{\ymin}{1e-6}
    \begin{subfigure}[b]{0.5\textwidth}
        \begin{tikzpicture}
\begin{axis}[%
width=0.951\fwidth,
height=\fheight,
at={(0\fwidth,0\fheight)},
scale only axis,
unbounded coords=jump,
title = error in expectation $e_E$,
xmin=\xmin,
xmax=10,
ymode=log,
ymin=\ymin,
ymax=\ymax,
legend style={
fill=none,
draw=none, 
text width=2.3cm, 
yshift = -0.25cm,
align=left,
font = \small,
},
legend pos = south west, 
xlabel = ROM dimension $r$,
yminorticks=true,
axis background/.style={fill=white},
ymajorgrids=true, 
xmajorgrids=true
]
    \addplot[mark=square*, cPOD] table [y=err1, x=rank] {data_edit2/Heat2d/E_error_POD.txt};     
    \addplot[mark=*, cOI1] table [y=err1, x=rank] {data_edit2/Heat2d/E_error_OpInf.txt}; 
    \addplot[mark=*, cOI2] table [y=err2, x=rank] {data_edit2/Heat2d/E_error_OpInf.txt};
    \addplot[mark=*, cOI3] table [y=err3, x=rank] {data_edit2/Heat2d/E_error_OpInf.txt};
    \legend{POD,OpInf $L=10^1$,OpInf $L=10^2$, OpInf $L=10^4$}
\end{axis}
\end{tikzpicture}%
    \end{subfigure}%
    \begin{subfigure}[b]{0.5\textwidth}
        \begin{tikzpicture}
\begin{axis}[%
width=0.951\fwidth,
height=\fheight,
at={(0\fwidth,0\fheight)},
scale only axis,
unbounded coords=jump,
title = error in covariance $e_C$,
xmin=\xmin,
xmax=\xmax,
ymode=log,
ymin=\ymin,
ymax=\ymax,
yminorticks=true,
xlabel = ROM dimension $r$,
axis background/.style={fill=white},
ymajorgrids=true, 
xmajorgrids=true
]
    \addplot[mark=square*, cPOD] table [y=err1, x=rank] {data_edit2/Heat2d/C_error_POD.txt};     
    \addplot[mark=*, cOI1] table [y=err1, x=rank] {data_edit2/Heat2d/C_error_OpInf.txt}; 
    \addplot[mark=*, cOI2] table [y=err2, x=rank] {data_edit2/Heat2d/C_error_OpInf.txt};
    \addplot[mark=*, cOI3] table [y=err3, x=rank] {data_edit2/Heat2d/C_error_OpInf.txt};
\end{axis}
\end{tikzpicture}%
    \end{subfigure}%
    \hfill
    \vspace{0.2cm}
     \begin{subfigure}[b]{0.5\textwidth}
        \begin{tikzpicture}
\begin{axis}[%
width=0.951\fwidth,
height=\fheight,
at={(0\fwidth,0\fheight)},
scale only axis,
unbounded coords=jump,
title = weak error $e_{\Phi,1}(T)$,
xmin=\xmin,
xmax=\xmax,
ymode=log,
ymin=\ymin,
ymax=\ymax,
yminorticks=true,
xlabel = ROM dimension $r$,
axis background/.style={fill=white},
ymajorgrids=true, 
xmajorgrids=true
]
    \addplot[mark=square*, cPOD] table [y=err1, x=rank] {data_edit2/Heat2d/f1_error_POD.txt};     
    \addplot[mark=*, cOI1] table [y=err1, x=rank] {data_edit2/Heat2d/f1_error_OpInf.txt}; 
    \addplot[mark=*, cOI2] table [y=err2, x=rank] {data_edit2/Heat2d/f1_error_OpInf.txt};
    \addplot[mark=*, cOI3] table [y=err3, x=rank] {data_edit2/Heat2d/f1_error_OpInf.txt};
\end{axis}
\end{tikzpicture}%
    \end{subfigure}%
    \begin{subfigure}[b]{0.5\textwidth}
        \begin{tikzpicture}
\begin{axis}[%
width=0.951\fwidth,
height=\fheight,
at={(0\fwidth,0\fheight)},
scale only axis,
unbounded coords=jump,
title = weak error $e_{\Phi,2}(T)$,
xmin=\xmin,
xmax=\xmax,
ymode=log,
ymin=\ymin,
ymax=\ymax,
yminorticks=true,
xlabel = ROM dimension $r$,
axis background/.style={fill=white},
ymajorgrids=true, 
xmajorgrids=true
]
    \addplot[mark=square*, cPOD] table [y=err1, x=rank] {data_edit2/Heat2d/f2_error_POD.txt};     
    \addplot[mark=*, cOI1] table [y=err1, x=rank] {data_edit2/Heat2d/f2_error_OpInf.txt}; 
    \addplot[mark=*, cOI2] table [y=err2, x=rank] {data_edit2/Heat2d/f2_error_OpInf.txt};
    \addplot[mark=*, cOI3] table [y=err3, x=rank] {data_edit2/Heat2d/f2_error_OpInf.txt};
\end{axis}
\end{tikzpicture}%
    \end{subfigure}%
    \caption{$2d$ Heat equation. The subspace computation, data generation and testing were performed with $s=100$ time steps of size $h=10^{-2}$. The number of samples drawn to estimate the error quantities was $10^6$.}
    \label{fig:2dHeat}
\end{figure}

\section{Conclusion and further research}
The inclusion of even a ``simple'' noise in the system dynamics introduced several new challenges in the OpInf approach, making a direct application of the OpInf method, developed in \cite{peherstorfer_data-driven_2016}, is infeasible.
In this paper, a nonintrusive reduced order modelling method was developed that approximates the FOM in the weak sense, i.e., in distribution.
It was shown that the inferred coefficients can be made arbitrarily close to the intrusive POD coefficients by chosing appropriate $L,h$, and $r$.
Together with the continuity of the expectation and covariance of linear SDEs with additive Gaussian noise with regard to the initial condition and system coefficients this implies that the presented method provides a ROM which is close in distribution to the POD-ROM.
Regarding the subspaces spanned left singular vectors corresponding to non-zero singular values, it was shown that the subspace obtained from the moment-snapshots is always contained in the subspace obtained from the state snapshots.
Although the introduced method performed well in the numerical experiments of Section \ref{sec:Experiments}, several challenges remain to be addressed.
First, estimating empirical moments requires many more FOM trajectories than in the deterministic setting, making the approach computationally expensive given the high cost of obtaining each trajectory.
However, since the presented method achieved relatively good results already for $L=10$, one possible intermediate step could be to investigate the confidence bounds of the approximation errors with respect to the number of training samples.
Second, we see potential in incorporating knowledge about the FOM into the inference. 
If, for example the noise dimension $d$ is small, this could be achieved by enforcing a low rank condition in the covariance estimation and the diffusion inference.
Finally, while this inference method can technically be employed on SDEs with a polynomial drift of degree larger than 1, the resulting ROM is still only informed by the first two moments, even though the FOM state variable is not necessarily Gaussian anymore. 
Consequently, the OpInf ROM may approximate the FOM expectation and covariance well but fail to capture higher moments. 
Extending the framework to handle higher-order polynomial drifts or even multiplicative noise is therefore an interesting direction for future research.

\section*{Reproducibility}
The code used in this publication is available at \url{https://github.com/JMNicolaus/OperatorInference_for_SDEs}

\section*{Acknowledgments}
The research has been partially funded by the Deutsche Forschungsgemeinschaft (DFG) - Project-ID 318763901 - SFB1294 as well as by the DFG individual grant ``Low-order approximations for large-scale problems arising in the context of high-dimensional
PDEs and spatially discretized SPDEs''-- project number 499366908.
In particular, we would like to acknowledge the productive discussions at the annual SFB1294 Spring School 2023.
We especially thank Dr. Thomas Mach for the many insightful discussions.

\appendix
\section{A Gronwall lemma for integrable functions}
Various generalisations of the Gronwall lemma have been proposed \cite{willettLinearGeneralizationGronwall1965,willettNonlinearVectorIntegral1964,pachpatte_inequalities_1998}.
However, most sources only provide statements for continuous coefficients and state the extendability of the results to integrable functions by referencing \cite{beesack1975gronwall}, which is, to our knowledge, unavailable in digital format and not readily accessible. 
We therefore provide a version of the Gronwall lemma for integrable functions.
The proof follows the standard arguments for proving inequalities of this type. 
\begin{lemma}\label{lem:Gronwall}
Let $\alpha,\beta:[0,T]\to\R_+$ be real non-negative integrable functions. Furthermore, let $\alpha$ be non-decreasing.
Let $u\in L^1([0,T])$ be such that 
\begin{align*}
    u(t) \leq \alpha(t) + \int_0^t\beta(s)u(s)\rmd s,
\end{align*}
holds for all $t\in[0,T]$.
If $\beta u\in L^1([0,T])$ and $\beta\alpha \in L^1([0,T])$, 
then $u$ admits to the upper bound
\begin{align*}
    u(t) \leq  \alpha(t)\exp\left( \int_0^t \beta(s) \rmd s\right)
\end{align*}
for all $t\in[0,T]$.
\end{lemma}
\begin{proof}
    The proof makes extensive use of the fundamental theorem of analysis for the Lebesgue integral \cite{cohnMeasureTheorySecond2013}.
    First, since $\beta\in L^1([0,T])$, the function $t\mapsto\int_0^t \beta(s)\rmd s$ is absolutely continuous, almost everywhere differentiable and 
    \[\beta(t) = \left(\int_0^t\beta(s)\rmd s\right)'\]
    holds almost everywhere.
    If one defines 
    \begin{align*}
    v(t) =\exp\left( -\int_0^t \beta(s)\rmd s\right),    
    \end{align*} 
    then $v'(t) = -\beta(t)v(t)$ almost everywhere.
    By the same argument, the integral of $\beta u$ is absolutely continuous and differentiable almost everywhere. 
    Now, let $h(t) := v(t)\int_0^t\beta(s)u(s)\rmd s$. 
    Then, by applying the chain rule, it follows that
    \begin{align*}
         h'(t) &= -\beta(t)v(t)\int_0^t\beta(s)u(s)\rmd s + v(t)\beta(t)u(t)\\
        &=v(t)\beta(t)\left(u(t)-\int_0^t\beta(s)u(s)\rmd s\right)\\
        &\leq v(t)\beta(t)\alpha(t)
    \end{align*}
    for almost all $t\in[0,T]$ exploiting the assumption on $u$.
    Since $v$ is a continuous function and $\beta\alpha\in L^1([0,T])$,  their product $v\beta\alpha$ is integrable as well. 
    Hence, by the above inequality, it holds that 
    \begin{align*}
        h(t) = h(t)-h(0) = \int_0^t h'(s)\rmd s \leq \int_0^t v(s)\beta(s)\alpha(s)\rmd s
    \end{align*}
 for all $t\in [0, T]$. 
 Dividing by $v(t)$ and adding $\alpha(t)$ on both sides leads to 
    \begin{align*}
        u(t) \leq \alpha(t) + \int_0^t\beta(s)u(s)\rmd s \leq \alpha(t) + \int_0^t\beta(s)\alpha(s)\frac{v(s)}{v(t)}\rmd s 
    \end{align*}
for all $t\in [0, T]$. 
Since $\alpha$ is non-decreasing, the right-hand side can be bounded by
    \begin{align*}
        \alpha(t) + \int_0^t\beta(s)\alpha(s)\frac{v(s)}{v(t)}\rmd s\leq \alpha(t)\exp\left(\int_0^t\beta(s)\rmd s\right)\quad 
    \end{align*}
    for all $t\in[0,T]$, completing the proof.

\end{proof}

\bibliographystyle{plain}
\bibliography{references_without_url}

\end{document}